# TRACES OF SYMMETRIC MARKOV PROCESSES AND THEIR CHARACTERIZATIONS

By Zhen-Qing Chen,[1] Masatoshi Fukushima[2] and Jiangang Ying[3]

*University of Washington, Kansai University and Fudan University*

*Dedicated to the Memory of Martin L. Silverstein*

Time change is one of the most basic and very useful transformations for Markov processes. The time changed process can also be regarded as the trace of the original process on the support of the Revuz measure used in the time change. In this paper we give a complete characterization of time changed processes of an arbitrary symmetric Markov process, in terms of the Beurling–Deny decomposition of their associated Dirichlet forms and of Feller measures of the process. In particular, we determine the jumping and killing measure (or, equivalently, the Lévy system) for the time-changed process. We further discuss when the trace Dirichlet form for the time changed process can be characterized as the space of finite Douglas integrals defined by Feller measures. Finally, we give a probabilistic characterization of Feller measures in terms of the excursions of the base process.

## Contents

1. Introduction
2. Feller measures and trace of Dirichlet forms
3. Space of functions with finite Douglas integrals
4. Excursions and Feller measures

Received August 2004; revised March 2005.
[1]Supported in part by NSF Grant DMS-03-03310.
[2]Supported in part by Grant-in-Aid for Scientific Research of MEXT No. 15540142.
[3]Supported in part by NSFC No. 10271109.

Appendix: Lévy system and Beurling–Deny decomposition for symmetric right processes

**1. Introduction.** Time change is one of the most basic and very useful transformations for Markov processes, which has been studied by many authors. However, a precise characterization of the time-changed process of a symmetric Markov process $X$ on a state space $E$ by a Revuz measure whose quasi-support $F$ is a proper subset of $E$ has only been started very recently. In [19], Fukushima, He and Ying derived a characterization for the time-changed process of $X$ when $X$ is a conservative $m$-symmetric diffusion, $m(E) < \infty$ and $F$ is a closed set which is negligible by the energy measure of $X$. The time-changed process has $F$ as its state space so it can be regarded as the trace process of $X$ on $F$.

The following is a prototype of the problem that we will study in this paper. Suppose $X$ is a Lévy process in $\mathbb{R}^n$ that is the sum of a Brownian motion in $\mathbb{R}^n$ and an independent spherically symmetric $\alpha$-symmetric stable process in $\mathbb{R}^n$, where $n \geq 1$ and $\alpha \in (0,2)$. Denote by $B(x,r)$ the open ball in $\mathbb{R}^n$ centered at $x \in \mathbb{R}^n$ with radius $r$. Its Euclidean closure is denoted by $\overline{B(x,r)}$. Let $F = \overline{B(0,1)} \cup \partial B(x_0,1)$, where $x_0 \in \mathbb{R}^n$ with $|x_0| = 3$. What is the trace process of $X$ on the closed set $F$? More precisely, let $\mu(dx) := \mathbb{1}_{\overline{B(0,1)}}(x)\,dx + \sigma_{\partial B(x_0,1)}$, where $\sigma_{\partial B(x_0,1)}$ denotes the Lebesgue surface measure of $\partial B(x_0,1)$. It is easy to see that $\mu$ is a smooth measure of $X$ and it uniquely determines a positive continuous additive functional $A^\mu = \{A^\mu_t, t \geq 0\}$ of $X$ having $\mu$ as its Revuz measure. Define its inverse

$$\tau_t := \inf\{s > 0 : A^\mu_s > t\} \qquad \text{for } t \geq 0.$$

Then the time changed process $Y_t := X_{\tau_t}$ is a symmetric Markov process on $F$, which can be regarded as the trace process of $X$ on $F$. So the more precise question is the following:

*Question*: Can we characterize the time changed process $Y$?

As a special case of the main results obtained in this paper, we are able to answer this question by determining its Dirichlet form on $L^2(F,\mu)$. It is known that the Dirichlet form $(\mathcal{E},\mathcal{F})$ for $X$ on $L^2(\mathbb{R}^n,dx)$ is given by

$$\mathcal{F} = W^{1,2}(\mathbb{R}^n, dx) := \{u \in L^2(\mathbb{R}^n, dx) : \nabla u \in L^2(\mathbb{R}^n, dx)\}$$

$$\mathcal{E}(u,u) = \frac{1}{2}\int_{\mathbb{R}^n} |\nabla u(x)|^2\,dx + \mathcal{A}(n,-\alpha)\int_{\mathbb{R}^n \times \mathbb{R}^n} \frac{(u(x)-u(y))^2}{|x-y|^{n+\alpha}}\,dx\,dy$$

for $u \in \mathcal{F}$,

where

(1.1) $$\mathcal{A}(n,-\alpha) = \frac{|\alpha|2^{\alpha-1}\Gamma((\alpha+n)/2)}{\pi^{n/2}\Gamma(1-\alpha/2)}$$

is a positive constant that depends only on $n$ and $\alpha$. Let $(\check{\mathcal{E}}, \check{\mathcal{F}})$ be the symmetric Dirichlet form of $Y$ on $L^2(F, \mu)$. The following explicit Beurling–Deny decomposition for $Y$ follows directly from Theorems 2.7, 2.10 and 2.11 below. Let $(\mathcal{E}, \mathcal{F}_e)$ denote the extended Dirichlet space for $(\mathcal{E}, \mathcal{F})$. Then

$$(1.2) \qquad \check{\mathcal{F}} = \mathcal{F}_e|_F \cap L^2(F, \mu),$$

$$\check{\mathcal{E}}(u, u) = \frac{1}{2} \int_F |\nabla u(x)|^2 \, dx$$

$$(1.3) \qquad + \int_{F \times F} (u(x) - u(y))^2 \left( \frac{1}{2} U(dx, dy) + \frac{\mathcal{A}(n, -\alpha)}{|x-y|^{n+\alpha}} \, dx \, dy \right)$$

$$+ \int_F u(x)^2 V(dx)$$

for $u \in \check{\mathcal{F}}$. Here $U$ and $V$ are the Feller measure on $F \times F$ and the supplementary Feller measure on $F$, respectively, defined through the energy functional $L$ for the subprocess $X^0$ of $X$ killed upon leaving the open set $\mathbb{R}^n \setminus F$:

$$(1.4) \qquad \int_{F \times F} f(x)g(y) U(dx, dy) := L(\mathbf{H}f, \mathbf{H}g)$$

$$:= \lim_{t \downarrow 0} \frac{1}{t} \int_{\mathbb{R}^n \setminus F} (\mathbf{H}f(x) - P_t^0 \mathbf{H}f(x)) \mathbf{H}g(x) \, dx,$$

$$\int_F f(x) V(dx) := L(\mathbf{H}f, q)$$

$$(1.5) \qquad = \lim_{t \downarrow 0} \frac{1}{t} \int_{\mathbb{R}^n \setminus F} (\mathbf{H}f(x) - P_t^0 \mathbf{H}f(x)) q(x) \, dx,$$

$$q(x) = \mathbb{1} - \mathbf{H}\mathbb{1}(x)$$

for any nonnegative bounded Borel measurable functions $f$ and $g$ on $F$. Here $\{P_t^0, t > 0\}$ is the transition semigroup of the subprocess $X^0$ and

$$\mathbf{H}f(x) = \mathbf{E}_x[f(X_{\sigma_F})] \qquad \text{for } x \in \mathbb{R}^n \setminus F,$$

is the "harmonic" extension of $f$ in the open set $\mathbb{R}^n \setminus F$, where $\sigma_F := \inf\{t \geq 0 : X_t \in F\}$. The measure $U$ is named after W. Feller, who introduced such measure for discrete Markov chains when studying their boundary theory.

In fact, in this paper we will study the above type of the problem for a general irreducible $m$-symmetric Markov process $X$ on a general state space $E$ which not only can have discontinuous sample paths but also can have killings inside $E$ or have finite lifetime, and for any *quasi-closed* subset $F$ of $E$. It is important for $F$ being quasi-closed rather than closed, since the notion of being quasi-closed is invariant under quasi-homeomorphism, while the notion of being closed is not. (Quasi-homeomorphism, whose definition can be found in [8] or in the proof of Theorem A.1 in the Appendix, is an equivalent relation not only at the quadratic forms level but also at the processes level. So fundament analytic and probabilistic properties of symmetric Markov processes should be preserved or invariant under quasi-homeomorphisms.) Since $F$ is only assumed to be quasi-closed in this paper, we can address both aspects of the problem and can apply results from regular Dirichlet forms through quasi-homeomorphisms (cf. [8]). If we consider the time-change problem for the process $X$, the smooth measure $\mu$ for time change is given. In this case, we take $F$ to be the quasi-support of $\mu$ which is quasi-closed. If we consider the trace problem for the process $X$, a quasi-closed set $F$ is given. We point out that any nontrivial quasi-closed subset $F$ is the quasi-support of a smooth measure $\mu$ of $X$ and we can take one of these measures $\mu$ to do time change. In this case, we fix a smooth measure $\mu$ of $X$ whose quasi-support is $F$. We are able to give a complete characterization of the Dirichlet form for the time-changed process $Y$ obtained from $X$ through $\mu$ and derive its Beurling–Deny decomposition of the time changed process similar to the one given in (1.2) and (1.3) (see Theorems 2.7, 2.10 and 2.11 below).

In particular, we determine the jumping and killing measure (or, equivalently, the Lévy system) for the time-changed process. More precisely, let $J$ and $\kappa$ be the jumping measure and the killing measure of $X$. We show in Section 2 below that

$$(1.6) \qquad \check{J} := \tfrac{1}{2} U + J|_{F \times F} \quad \text{and} \quad \check{\kappa} := V + \kappa|_F$$

are the jumping measure and the killing measure of the time-changed process $Y$ of $X$. Note that the jumping measure $\check{J}$ and the killing measure $\check{\kappa}$ depend only on $X$ and the quasi-closed set $F$; they do not depend on the selection of the measure $\mu$ that has $F$ as its quasi-support. In fact, if $\mu_1$ and $\mu_2$ are two smooth measures both having $F$ as their quasi-support, then by a result due to Silverstein and Fitzsimmons (see [13] and [31]), the time-changed process of $X$ by the smooth measure $\mu_1$ is a time-change of the time-changed process of $X$ by the smooth measure $\mu_2$.

Our results extend the recent work by Fukushima, He and Ying [19] where only a conservative symmetric diffusion process $X$ is considered and $F$ is assumed to be closed. When $X$ is a conservative diffusion (where $J = \kappa = 0$) and under an additional condition for $X$ which is fulfilled when $m(E) < \infty$, (1.6) has been proved in [19], Chapter 5, by relating the jumping and killing measures of $Y$ to expectations of certain homogeneous random measures involving end points of excursions of $X$ away from the closed set $F$. Under the further condition that $F$ is negligible for the energy measure of $X$, it is shown in [19], Chapter 6, that the Dirichlet form for the time-changed process $Y$ has no strongly local part and, therefore, the Beurling–Deny decomposition just has the jumping part.

In this paper, we do not use excursions to derive (1.6). Instead, by sharpening and extending those computations performed in [19], Chapter 6, to the present general situation, we make a direct and precise analysis of the value $\mathcal{E}(\mathbf{H}u, \mathbf{H}u)$ for $u \in \check{\mathcal{F}}$, decomposing it into a sum of terms involving the strongly local part of the energy measure of $\mathbf{H}u$, the measures $J, \kappa$ and the Feller measures $U, V$, eventually leading us to the Beurling–Deny decomposition of the trace Dirichlet form $\check{\mathcal{E}}$. One of the key steps is to identify the strongly local part $\check{\mathcal{E}}^c(u,u)$ in the Beurling–Deny decomposition of $\check{\mathcal{E}}(u,u)$ with $\mu^c_{\langle \mathbf{H}u \rangle}(F)$, where $\mu^c_{\langle \mathbf{H}u \rangle}$ is the strongly local part of the energy measure of $\mathbf{H}u$ with respect to $(\mathcal{E}, \mathcal{F})$. The stochastic calculus for general martingales with possibly discontinuous sample paths plays a key role in our approach.

When $X$ is the $n$-dimensional Brownian motion and $F$ is a compact hypersurface of class $C^3$, an explicit expression of the Feller measure $U$ is exhibited in [19], Example 2.1. We shall give an explicit expression of the supplementary Feller measure $V$ in this example at the end of Section 2.

The identification (1.6) particularly implies that the (generalized) Douglas integral

$$(1.7) \qquad \tfrac{1}{2} \int_{F \times F} (\varphi(\xi) - \varphi(\eta))^2 U(d\xi, d\eta) + \int_F \varphi(\xi)^2 V(d\xi)$$

is finite for any function $\varphi$ in the extended trace Dirichlet space $\check{\mathcal{F}}_e$. In Section 3, we are concerned with conditions to ensure that the space $\check{\mathcal{F}}_e$ coincides with the space of functions with finite Douglas integrals. The conditions will be given in relation to the reflected Dirichlet space $(\mathcal{E}^{\mathrm{ref}}, (\mathcal{F}^0)^{\mathrm{ref}})$ of the part $(\mathcal{E}^0, \mathcal{F}^0)$ of the Dirichlet form $(\mathcal{E}, \mathcal{F})$ on the set $E_0 = E \setminus F$. $(\mathcal{E}^0, \mathcal{F}^0)$ is associated with the absorbed process $X^0$ obtained from $X$ by killing upon leaving the set $E_0$.

The notion of the reflected Dirichlet space was introduced by Silverstein (cf. [31]) and further studied by the first author [6] for a general transient regular Dirichlet space, while the space of functions of finite Douglas integrals was studied by Doob and the second author [17] in a different setting for the absorbed Brownian motion $X^0$ on an Euclidean domain $E_0$ and its Martin boundary $F$.

Since $E_0$ is only quasi-open, the Dirichlet form $(\mathcal{E}^0, \mathcal{F}^0)$ on $L^2(E_0; m)$ is no longer regular in general, but it is quasi-regular, as will be shown in Lemma 2.2. Based on this fact, we are able to extend the definition of the reflected Dirichlet space given in [6] to $(\mathcal{E}^0, \mathcal{F}^0)$ by making use of the notion of the energy functional $L$ of excessive functions for the process $X^0$.

As applications, we present examples of the trace Dirichlet spaces for the reflecting Brownian motion $X$ and the reflected symmetric stable process $X$ on a Euclidean domain at the end of Section 3.

By reversing the argument in [19], we may well derive from the identification (1.6) expressions of the Feller measures $U, V$ in terms of homogeneous random measures generated by end points of excursions of $X$ away from $F$. This is what will be done in Section 4. A direct derivation of such expressions seems to be hard in the present generality unless $X$ is conservative as in the case of [19].

In the Appendix we show that a Lévy system exists for any symmetric right process associated with a quasi-regular Dirichlet form $(\mathcal{E}, \mathcal{F})$ and that the probabilistic characterization of the Beurling–Deny decomposition (including the jumping measure and killing measure) of $(\mathcal{E}, \mathcal{F})$ remain true for quasi-regular Dirichlet forms. This result might be known to the experts, but we are unable to find an exact reference for it. Since this result is used in the paper, for the reader's convenience, we record it in the Appendix. We further show that this probabilistic characterization is independent of the choice of a particular process associated with the Dirichlet form. These results are used in Section 4. Some basic $\mathcal{E}$-quasi notions, such as $\mathcal{E}$-nest, $\mathcal{E}$-polar, $\mathcal{E}$-quasi-everywhere, $\mathcal{E}$-quasi-closed set, $\mathcal{E}$-quasi-continuous, quasi-homeomorphism between Dirichlet forms, etc., are also reviewed in this Appendix.

We point out that under some extra condition, LeJan has obtained in Section 3 of [27] the same results as Section 2 of the present paper for a Hunt process $X$ associated with a nonsymmetric sectorial regular Dirichlet form and for a closed set $F$. Along with [26], nice potential theoretic methods were systematically utilized in [27] under the condition that the Dirichlet space is continuously embedded into $L^2(E; m)$. This condition, however, excludes many interesting examples, such as reflecting Brownian motion in a unit disk while $F$ is the unit circle. In [24] Kunita treated Douglas integral representations on the Martin boundary for multidimensional diffusions in the same spirit as Doob [11] and Fukushima [17]. For other approaches, such as those using the framework of balayage spaces and harmonic spaces, to the problems related to the traces of Markov processes and their potential theory, see [3, 23] and the references therein.

In this paper we use ":=" as a way of definition, which is read as "is defined to be." For a real-valued function $\varphi(t)$ on $\mathbb{R}$, $\varphi(t-) := \lim_{s \uparrow t} \varphi(s)$ denotes its left-hand limit at $t$ if it exists. For two real numbers $a$ and $b$, $a \vee b := \max\{a, b\}$ and $a \wedge b := \min\{a, b\}$. For a locally compact metric space $E$, $C_c(E)$ denotes the space of real-valued continuous functions on $E$ with compact support. For a Borel subset $K$ of $E$, we will use $\mathcal{B}(K)$, $\mathcal{B}(K)^+$ and $\mathcal{B}(K)_b^+$ to denote the space of Borel measurable functions on $K$, the space of nonnegative Borel measurable functions on $K$ and the space of nonnegative bounded Borel measurable functions on $K$, respectively.

**2. Feller measures and trace of Dirichlet forms.** Throughout this paper, let $(\mathcal{E}, \mathcal{F})$ be an irreducible quasi-regular symmetric Dirichlet form on $L^2(E,m)$, where $E$ is a Hausdorff metric space and the measure $m$ has $\mathrm{supp}[m] = E$. Let $X$ be the $m$-symmetric right process associated with $(\mathcal{E}, \mathcal{F})$, whose life time will be denoted as $\zeta$.

In view of the quasi-homeomorphism method in [8], without loss of generality, we may and do assume that $E$ is a locally compact separable metric space, $m$ is a positive Radon measure on $E$ with $\mathrm{supp}[m] = E$, $(\mathcal{E}, \mathcal{F})$ is an irreducible regular symmetric Dirichlet form in $L^2(E,m)$, and $X = (X_t, \mathbf{P}_x)$ is an $m$-symmetric Hunt process associated with $(\mathcal{E}, \mathcal{F})$. We will use $(\mathcal{E}, \mathcal{F}_e)$ to denote the extended Dirichlet space of $(\mathcal{E}, \mathcal{F})$ and $\mathcal{E}_1 := \mathcal{E} + (\cdot, \cdot)_{L^2(E,m)}$. The expectation with respect to the probability measure $\mathbf{P}_x$ will be denoted as $\mathbf{E}_x$. Throughout this paper, we use the convention that any function takes value $0$ at the cemetery point $\partial$ added to $E$. For basic notions in Dirichlet forms, such as nest, capacity, quasi-everywhere (abbreviated q.e.), quasi-continuous, etc., we refer the reader to [20] and [28]. (See also the Appendix of this paper.) Every element $u$ in $\mathcal{F}_e$ admits a quasi-continuous version. We assume throughout this paper that functions in $\mathcal{F}_e$ are always represented by their quasi-continuous version. In the sequel, the abbreviations CAF, PCAF and MAF stands for "continuous additive functional," "positive continuous additive functional" and "martingale additive functional," respectively, whose definitions can be found in [20].

We prove a lemma that will be used later.

LEMMA 2.1. *Let $\mu$ be a smooth measure with $\mu(E) < \infty$ and let $A^\mu$ be the PCAF of $X$ with Revuz measure $\mu$. Then*

$$\lim_{t \downarrow 0} \frac{1}{t} \mathbf{E}_m[A_t^\mu; t \geq \zeta] = 0.$$

PROOF. Let $P_t f(x) := \mathbf{E}_x[f(X_t)]$. For $t > 0$, by the Markov property of $X$ and Theorem 5.1.3(iii) of [20],

$$\frac{1}{t} \mathbf{E}_m[A_t^\mu; t \geq \zeta]$$
$$= \frac{1}{t} \mathbf{E}_m[A_t^\mu] - \frac{1}{t} \mathbf{E}_m[A_t^\mu; t < \zeta]$$
$$= \frac{1}{t} \mathbf{E}_m[A_t^\mu] - \frac{1}{t} \mathbf{E}_m\left[\int_0^t \mathbb{1}_{\{t < \zeta\}} dA_s^\mu\right] ds$$
$$= \frac{1}{t} \mathbf{E}_m[A_t^\mu] - \frac{1}{t} \mathbf{E}_m\left[\int_0^t P_{t-s}\mathbb{1}(X_s) dA_s^\mu\right]$$
$$= \frac{1}{t} \mathbf{E}_m[A_t^\mu] - \frac{1}{t} \int_0^t \left(\int_E P_s\mathbb{1}(x) P_{t-s}\mathbb{1}(x) \mu(dx)\right) ds$$
$$= \frac{1}{t} \mathbf{E}_m[A_t^\mu] - \int_E \left(\frac{1}{t} \int_0^t P_{t-s}\mathbb{1}(x) P_s\mathbb{1}(x) ds\right) \mu(dx).$$

Since $\mu(E) < \infty$, by the dominated convergence theorem,

$$\lim_{t \downarrow 0} \frac{1}{t} \mathbf{E}_m[A_t^\mu; t \geq \zeta]$$
$$= \mu(E) - \int_E \left( \lim_{t \downarrow 0} \frac{1}{t} \int_0^t P_{t-s} \mathbb{1}(x) P_s \mathbb{1}(x) \, ds \right) \mu(dx)$$
$$= \mu(E) - \mu(E) = 0.$$

The lemma is proved. $\square$

For a closed set $F \subset E$, define

$$\mathcal{F}_F := \{ u \in \mathcal{F} : u = 0 \text{ $m$-a.e. on } E \setminus F \}.$$

Denote by Cap the $\mathcal{E}_1$-capacity defined by the form $(\mathcal{E}, \mathcal{F})$. The terms such as quasi-continuous functions, quasi-everywhere (q.e. in abbreviation), quasi-closed sets and generalized nests will be used exclusively in relation to this capacity. It is known (e.g., [18], Lemma 2.1, [10], Lemma 2.1) that those classical notions for the regular Dirichlet form $(\mathcal{E}, \mathcal{F})$ can be identified with the $\mathcal{E}$-quasi-notions of [28] (see Appendix for their definitions) as follows: an increasing sequence of closed subsets of $E$ is a generalized nest if and only if it is an $\mathcal{E}$-nest; a subset of $E$ is of zero capacity if and only if it is $\mathcal{E}$-exceptional; a numerical function defined q.e. on $E$ is quasi continuous if and only if it is $\mathcal{E}$-quasi continuous. Moreover, a set $F \subset E$ is quasi-closed if and only if it is $\mathcal{E}$-quasi-closed in the following sense (see also the Appendix below). A subset $F$ of $E$ is $\mathcal{E}$-quasi-closed if and only if there is an increasing sequence $\{K_n, n \geq 1\}$ of compact subsets of $E$ such that $\bigcup_{n \geq 1} \mathcal{F}_{K_n}$ is dense in $(\mathcal{F}, \mathcal{E}_1)$ and $F \cap K_n$ is closed for every $n \geq 1$.

Let $F$ be a quasi-closed subset $F$ of $E$ such that

$$\text{Cap}(F) > 0. \tag{2.1}$$

The notion of being quasi-closed is invariant under the quasi-homeomorphism of Dirichlet forms (see [8]). Furthermore, a quasi-closed set is q.e. finely closed in the sense that there is a properly exceptional set $N$ such that $F \setminus N$ is nearly Borel measurable and finely closed with respect to $X$ (cf. [20], Chapter 4.6). Since we are only concerned with assertions holding q.e., we may and do make a convention that the quasi-closed set $F$ is nearly Borel and finely closed already.

Let $E_0 = E \setminus F$. Under the present convention, $E_0$ is nearly Borel and finely open with respect to $X$. The subprocess of $X$ killed upon leaving $E_0$ will be denoted by $X^0$. To be more precise, we let

$$\tau_0 := \tau_{E_0} = \inf\{t \in [0, \zeta] : X_t \notin E_0\}, \tag{2.2}$$

so that

$$\tau_0 = \sigma_F \wedge \zeta, \qquad \mathbf{P}_x\text{-a.s. for } x \in E_0,$$

where $\sigma_F := \inf\{t > 0 : X_t \in F\}$. The subprocess $X^0$ is then defined by $X^0 = (X_t^0, \zeta^0, \mathbf{P}_x)_{x \in E_0}$, where

$$\zeta^0 =: \tau_0 \quad \text{and} \quad X_t^0 = \begin{cases} X_t & \text{for } t < \zeta^0, \\ \partial, & \text{for } t \geq \zeta^0. \end{cases}$$

The process $X^0$ is an $m$ symmetric standard process on $E_0$ and its Dirichlet form $(\mathcal{E}^0, \mathcal{F}^0)$ on $L^2(E_0, m)$ can be identified with the following space (see [20], Chapter 4.4):

(2.3) $\qquad \mathcal{F}^0 = \{u \in \mathcal{F} : u = 0 \text{ q.e. on } F\} \quad \text{and} \quad \mathcal{E}^0 = \mathcal{E}|_{\mathcal{F}^0 \times \mathcal{F}^0}.$

We use the terms $\mathcal{E}^0$-nest, $\mathcal{E}^0$-quasi continuous, $\mathcal{E}^0$-exceptional, and so on for those quasi notions exclusively related to the Dirichlet form $(\mathcal{E}^0, \mathcal{F}^0)$ in the sense of [8] or [28]. The resolvent of $X^0$ will be denoted by $G_\alpha^0$.

LEMMA 2.2. (i) *An increasing sequence $\{A_n\}$ of relatively closed subsets of $E_0$ is an $\mathcal{E}^0$-nest if and only if*

(2.4) $\qquad \mathbf{P}_x\left(\lim_{n \to \infty} \sigma_{E_0 \setminus A_n} < \tau_0\right) = 0 \qquad \text{for q.e. } x \in E_0.$

*Any $\mathcal{E}^0$-exceptional set $N \subset E_0$ is $X^0$-exceptional in the sense that there is a nearly Borel set $\widetilde{N}$ containing $N$ such that*

$$\mathbf{P}_x(\sigma_{\widetilde{N}} < \tau_0) = 0 \qquad \text{for } m\text{-a.e. } x \in E_0.$$

(ii) *If an increasing sequence of closed subsets $\{B_n\}$ of $E$ is a generalized nest (or, equivalently, an $\mathcal{E}$-nest), then $\{B_n \cap E_0\}$ is an $\mathcal{E}^0$-nest. The restriction to $E_0$ of a quasi-continuous function on $E$ is $\mathcal{E}^0$-quasi-continuous.*

(iii) *A set $N \subset E_0$ is $\mathcal{E}^0$-exceptional if and only if $\mathrm{Cap}(N) = 0$.*

(iv) *The Dirichlet form $(\mathcal{E}^0, \mathcal{F}^0)$ is transient and its extended Dirichlet space $\mathcal{F}_e^0$ can be described as*

(2.5) $\qquad \mathcal{F}_e^0 = \{u \in \mathcal{F}_e : u = 0 \text{ q.e. on } F\}.$

(v) *$X^0$ is a special standard process on $E_0$ and $(\mathcal{E}^0, \mathcal{F}^0)$ is a quasi-regular Dirichlet form on $L^2(E_0, m)$.*

PROOF. (i) Note that (2.4) is equivalent to

(2.6) $\qquad \sigma_{E \setminus A_n} \uparrow \sigma_F, \qquad \mathbf{P}_x\text{-a.s. for q.e. } x \in E_0.$

For a nearly Borel set $A \subset E$, let

$$\mathcal{F}_A = \{u \in \mathcal{F} : u = 0 \text{ q.e. on } E \setminus A\},$$

$$G_1^A f(x) = \mathbf{E}_x\left(\int_0^{\sigma_{E \setminus A}} e^{-t} f(X_t)\, dt\right), \qquad x \in E.$$

Then $G_1^A(\mathcal{B}^+(E) \cap L^2(E; m))$ is $\mathcal{E}_1$ dense in $\mathcal{F}_A$ by [20], Chapter 4.4. So (2.6) implies that $G_1^{A_n} f \in \mathcal{F}_{A_n}$ increases to $G_1^0 f$ $m$-a.e. and converges also in the $\mathcal{E}_1$-metric for $f \in \mathcal{B}^+(E) \cap L^2(E; m)$ and, accordingly,

(2.7) $$\bigcup_{n \geq 1} \mathcal{F}_{A_n} \text{ is } \mathcal{E}_1\text{-dense in } \mathcal{F}^0.$$

This proves that $\{A_n\}$ is an $\mathcal{E}^0$-nest.

Conversely, for $\sigma = \lim_{n \to \infty} \sigma_{E \setminus A_n}$ and $f \in \mathcal{B}^+(E) \cap L^2(E; m)$,

$$v(x) = \mathbf{E}_x\left(\int_\sigma^{\sigma_F} e^{-t} f(X_t)\, dt\right)$$

is a function in $\mathcal{F}^0$ that is $\mathcal{E}_1$-orthogonal to $\bigcup_{n=1}^\infty \mathcal{F}_{A_n}$ according to [20], Chapter 4.4. Hence, (2.7) implies $v = 0$ q.e. and (2.6) follows.

The second statement of (i) follows immediately from (2.4).

(ii) By [20], Lemma 5.1.6 or [28], $\{B_n\}$ is a generalized nest (or, equivalently, an $\mathcal{E}$-nest) if and only if

$$\mathbf{P}_x\left(\lim_{n \to \infty} \sigma_{E \setminus B_n} < \zeta\right) = 0 \qquad \text{for q.e. } x \in E.$$

So (2.4) holds for $A_n = B_n \cap E_0$. Since any increasing sequence $\{A_n\}$ of closed sets of $E$ with $\text{Cap}(E \setminus A_n) \downarrow 0$ is a generalized nest, the second statement of (ii) now follows.

(iii) The "if" part is another consequence of (ii). To prove the "only if" part, suppose a set $N \subset E_0$ is $\mathcal{E}^0$-exceptional but $\text{Cap}(N) > 0$. Then, in view of [20], Theorem 2.2.3, there exists a smooth measure $\mu$ on $E$ such that $\mu(N) > 0, \mu(E \setminus N) = 0$. Let $A_t$ be a PCAF of $X$ with Revuz measure $\mu$. Then $A_t$ is a PCAF in the strict sense of $X|_{E \setminus N_0}$ with Revuz measure $\mu|_{E \setminus N_0}$ for a properly exceptional set $N_0$.

Since $E_0$ is finely open, we can use [15], Theorem 2.22 to conclude that $t \mapsto A_{t \wedge \tau_0}$ is then a PCAF in the strict sense of the process $X^0|_{E_0 \setminus N_0}$ with Revuz measure $\mu|_{E_0 \setminus N_0}$, which charges no $X^0$-exceptional set accordingly. In particular, $\mu(N \setminus N_0) = 0$. Since $\mu(N_0) = 0$, we arrive at $\mu(N) = 0$, which is a contradiction. This proves that any $\mathcal{E}^0$-exceptional set $N$ satisfies $\text{Cap}(N) = 0$.

(iv) Since $(\mathcal{E}, \mathcal{F})$ is assumed to be irreducible, the transience of $(\mathcal{E}^0, \mathcal{F}^0)$ follows from the assumption (2.1) and [20], Theorem 4.6.6, Lemma 1.6.5.

The identification (2.5) of its extended Dirichlet space is shown in [20], Theorem 4.4.4, when $(\mathcal{E}, \mathcal{F})$ is transient. In general, we can show it by reducing the situation to the transient case. Here we show the inclusion $\supset$ in (2.5) since we shall use this inclusion only. Take any $u \in \mathcal{F}_e$ vanishing q.e. on $F$. On account of [20], Lemma 1.6.7, there exists a function $g$ satisfying condition (1.6.14) of [20] such that $u$ belongs to the extended Dirichlet space $\mathcal{F}_e^g$ of the perturbed Dirichlet form $(\mathcal{E}^g, \mathcal{F})$ which is transient. Therefore,

$$u \in (\mathcal{F}_e^g)_{E_0} = (\mathcal{F}_{E_0}^g)_e \subset \mathcal{F}_e^0,$$

as was to be proved.

(v) Since $F$ is finely closed, $X_{\sigma_F} \in F$. Hence, it follows from [4], Chapter IV, equation (4.33), that $X^0$ is special. The second assertion follows from [28], Theorem 5.1. □

The last statement (v) in the above lemma will be used in the next section.

Recall that, by our convention, every function in the space $\mathcal{F}^0$ and $\mathcal{F}_e^0$ is $\mathcal{E}^0$-quasi-continuous by Lemma 2.2(ii). We call a positive Borel measure $\mu$ on $E_0$ $\mathcal{E}^0$-*smooth* if $\mu$ charges no $\mathcal{E}^0$-exceptional set and there exists an $\mathcal{E}^0$-nest $\{A_n\}$ such that $\mu(A_n) < \infty$ for each $n$. Let $S^{(0)}$ denote the family of all $\mathcal{E}^0$-smooth measures on $E_0$. The restriction to $E_0$ of any smooth measure on $E$ is a member of $S^{(0)}$ in view of Lemma 2.2(ii)–(iii). A measure $\mu$ in $S^{(0)}$ is said to be *of finite 0-order energy integral* if there is a constant $C > 0$ such that

$$\int_{E_0} |v(x)|\mu(dx) \leq C\sqrt{\mathcal{E}(v,v)} \qquad \text{for every } v \in \mathcal{F}_e^0;$$

or, equivalently, there exists a function $G^0\mu \in \mathcal{F}^0$ such that

(2.8) $$\mathcal{E}(G^0\mu, v) = \int_{E_0} v(x)\mu(dx) \qquad \text{for every } v \in \mathcal{F}^0.$$

Such function $G^0\mu$ is then unique and called the 0-*order potential* of $\mu$. The totality of $\mathcal{E}^0$-smooth measures on $E_0$ of finite 0-order energy integrals will be denoted by $S_0^{(0)}$.

With these preparations, let us now consider the notion of the energy functional $L$ for the process $X^0$. The semigroup and the resolvent of $X^0$ will be denoted by $\{P_t^0, t \geq 0\}$ and $\{G_\alpha^0, \alpha > 0\}$, respectively. The 0-order potential for $X^0$ will be denoted as $G^0$. Denote by $\mathcal{S}_{E_0}$ the family of $m$-almost $X^0$-excessive functions finite $m$-a.e. on $E_0$. Let

(2.9) $$L(f,g) := \uparrow \lim_{t \downarrow 0} \frac{1}{t} \int_{E_0} (f - P_t^0 f)(x) g(x) m(dx)$$

be the energy functional of $f, g \in \mathcal{S}_{E_0}$. Here $\uparrow \lim_{t \downarrow 0}$ means that it is an increasing limit as $t \downarrow 0$. Equivalently,

(2.10) $$L(f,g) := \uparrow \lim_{\alpha \uparrow \infty} \alpha \int_{E_0} (f - \alpha G_\alpha^0 f)(x) g(x) m(dx),$$

where $\uparrow \lim_{\alpha \uparrow \infty}$ means that it is an increasing limit as $\alpha \uparrow \infty$. Note that $L(f,g)$ can be well defined without assuming the finiteness of $g$. We shall use this fact in the next section without an explicit mentioning. Note further that the 0-order potential $G^0 \nu$ of measure $\nu \in S_0^{(0)}$ is an element not only of $\mathcal{F}_e^0$ but also of $\mathcal{S}_{E_0}$, as we can easily see from [20], Theorem 2.2.1, and (2.8).

LEMMA 2.3. *Let $f, g \in \mathcal{S}_{E_0}$.*

(i) *If $\int_{E_0} f(x)g(x)m(dx) < \infty$, then $L(f,g) = L(g,f)$.*
(ii) *For any $\nu \in S_0^{(0)}$,*

$$\tag{2.11} L(f, G^0 \nu) = \int_{E_0} f(x) \nu(dx).$$

PROOF. (i) For $f, g \in \mathcal{S}_{E_0}$ with $fg \in L^1(E_0, m)$, we have, for each $t > 0$,

$$\frac{1}{t} \int_{E_0} (f - P_t^0 f)(x) g(x) m(dx) = \frac{1}{t} \int_{E_0} f(x) (g - P_t^0 g)(x) m(dx).$$

Passing $t \to 0$ yields $L(f, g) = L(g, f)$.

(ii) By the transience of $(\mathcal{E}, \mathcal{F}^0)$, there is a bounded, $L^1(E_0, m)$-integrable function $h$ that is strictly positive $m$-a.e on $E$ and $\int_{E_0} h(x) G^0 h(x) m(dx) < \infty$. Put $f_n = f \wedge (nh)$. Then $G^0 f_n, G^0 G_\alpha^0 f_n \in \mathcal{F}_e^0$ and

$$(f_n - \alpha G_\alpha^0 f_n, G^0 \nu)_{L^2(E_0, m)} = \mathcal{E}(G^0 f_n, G^0 \nu) - \alpha \mathcal{E}(G^0 G_\alpha^0 f_n, G^0 \nu)$$

$$= \int_{E_0} (G^0 f_n - \alpha G^0 G_\alpha^0 f_n)(x) \nu(dx)$$

$$= \int_{E_0} G_\alpha^0 f_n(x) \nu(dx).$$

By letting $n \to \infty$, we get

$$\tag{2.12} (f - \alpha G_\alpha^0 f, G^0 \nu)_{L^2(E_0, m)} = \int_{E_0} G_\alpha^0 g(x) \nu(dx).$$

Next, for $\alpha > 0$, put

$$f_\alpha = g - \alpha G_\alpha^0 f, \qquad f_{\alpha,n} = f_\alpha \wedge (nh).$$

Then

$$\int_{E_0} G^0 f_{\alpha,n} \nu(dx) = \mathcal{E}(G^0 f_{\alpha,n}, G^0 \nu) = (f_{\alpha,n}, G^0 \nu)_{L^2(E_0, m)}.$$

Passing $n \to \infty$, we get

$$\int_{E_0} G^0 f_\alpha(x) \nu(dx) = (f_\alpha, G^0 \nu)_{L^2(E_0, m)}.$$

Hence, by (2.12),

$$(2.13) \qquad (f - \alpha G_\alpha^0 f, G^0 \nu)_{L^2(E_0, m)} = (f_\alpha, G^0 \nu)_{L^2(E_0, m)}.$$

Multiplying both sides of (2.13) by $\alpha$ and then letting $\alpha \to \infty$ establishes (2.11) for the measure $\nu$ that is of finite energy with respect to $(\mathcal{E}^0, \mathcal{F}^0)$. □

For $\alpha \geq 0$, let $\mathbf{H}^\alpha$ denote the $\alpha$-order hitting measure of $F$; that is,

$$\mathbf{H}^\alpha(x, B) = \mathbf{E}_x[e^{-\alpha \tau_0} \mathbb{1}_B(X_{\tau_0}); \tau_0 < \infty] \qquad \text{for } x \in E_0 \text{ and } B \in \mathcal{B}(F).$$

When $\alpha = 0$, we will use $\mathbf{H}$ to denote $\mathbf{H}^0$. Since $F$ is a finely closed set, $\mathbf{H}^\alpha(x, \cdot)$ is carried by $F$. For $f \in \mathcal{B}(F)^+$, define

$$\mathbf{H}^\alpha f(x) := \mathbf{E}_x[e^{-\alpha \tau_0} f(X_{\tau_0}); \tau_0 < \infty] \qquad \text{for } x \in E.$$

Clearly, for any $f \in \mathcal{B}(F)^+$ and $\alpha \geq 0$, $\mathbf{H}^\alpha f$ is $\alpha$-excessive with respect to the subprocess $X^0$. Moreover, for every $f \in \mathcal{B}(F)_b^+$ and every $x \in E_0$,

$$(2.14) \qquad \lim_{t \to \infty} P_t^0 \mathbf{H} f(x) = \lim_{t \to \infty} \mathbf{E}_x[f(X_{\tau_0}); t < \tau_0 < \infty] = 0.$$

By Theorem 4.6.5 of [20], for $u \in \mathcal{F}_e$, $\mathbf{H}|u| < \infty$ q.e. on $E$ and $\mathbf{H}u \in \mathcal{F}_e$. Recall that, for $u \in \mathcal{F}_e$, the following Fukushima's decomposition holds uniquely (cf. [20]):

$$u(X_t) - u(X_0) = M^u + N^u \qquad \text{for } t \geq 0,$$

where $M^u$ is a martingale additive functional of $X$ having finite energy and $N^u$ is a continuous additive functional of $X$ having zero energy. In the sequel, we will use $\mu_{\langle u \rangle}$ to denote the Revuz measure for the predictable quadratic variation $\langle M^u \rangle$ for the square integrable martingale $M^u$.

Let $(N(x, dy), H)$ denote a Lévy system for the $m$-symmetric Hunt process $X$ on $E$. (See the Appendix below for the definition and properties of a Lévy system for symmetric Borel right processes.) The Revuz measure of the PCAF $H$ of $X$ will be denoted as $\mu_H$. We define

$$J(dx, dy) = \tfrac{1}{2} N(x, dy) \mu_H(dx) \quad \text{and} \quad \kappa(dx) = N(x, \partial) \mu_H(dx)$$

as the jumping measure and the killing measure of $X$ [or, equivalently, of $(\mathcal{E}, \mathcal{F})$]. The following Lemma 2.4 is an extension of Lemma 6.1 in [19]. Note here we do not assume $X$ being conservative.

Define $q(x) := 1 - \mathbf{H}\mathbb{1}(x) = \mathbf{P}_x(\tau_0 \geq \zeta)$. For $f, g \in \mathcal{B}(F)_b^+$, define

$$(2.15) \qquad U(f \otimes g) := L(\mathbf{H}f, \mathbf{H}g) \quad \text{and} \quad V(f) := L(\mathbf{H}f, q).$$

By Lemma 2.3(i), $U$ is a symmetric measure on $F \times F$, which will be called the *Feller measure* for $F$. The measure $V$ on $F$ will be called the *supplementary Feller measure* for $F$.

We will first assume the following condition:

(2.16) $$m(E_0) < \infty.$$

We will show in Theorem 2.11 below that this condition (2.16) can be dropped.

LEMMA 2.4. *Assume condition* (2.16) *holds. For any* $u \in \mathcal{F}_{e,b}$, *let* $w = \mathbf{H}(u^2) - (\mathbf{H}u)^2$. *Then* $w \in \mathcal{F}_{e,b}^0$ *and* $w = G^0\nu$ *with* $\nu = \mu_{\langle \mathbf{H}u \rangle}|_{E_0} \in S_0^{(0)}$. *Furthermore,*

(2.17) $$\mu_{\langle \mathbf{H}u \rangle}(E_0) + \int_{E_0} (\mathbf{H}u)^2(x)\kappa(dx)$$
$$= \lim_{\alpha \to \infty} \alpha(\mathbf{H}^\alpha \mathbb{1}, w)_{L^2(E_0, m)} + \int_F u(x)^2 V(dx).$$

PROOF. The proof of the first assertion is the same as that for Lemma 6.1 in [19], although we now work with the space $S_0^{(0)}$ of measures on the quasi-open set $E_0$ rather than an open set. So it is omitted here. Note that we only need the inclusion $\supset$ in relation (2.5) in our proof.

Since $w$ is bounded and $m(E_0) < \infty$, we have, by Lemma 2.3,

(2.18) $$\nu(E_0) = L(1, w) = L(w, 1)$$
$$= L(w, \mathbf{H}\mathbb{1}) + L(w, q) = L(\mathbf{H}\mathbb{1}, w) + L(w, q)$$

and, hence,

(2.19) $$\mu_{\langle \mathbf{H}u \rangle}(E_0) = \lim_{\alpha \to \infty} \alpha(\mathbf{H}^\alpha \mathbb{1}, w)_{L^2(E_0, m)} + L(w, q).$$

On the other hand,

(2.20) $$L(w, q) = L(\mathbf{H}(u^2), q) - \lim_{t \downarrow 0} \frac{1}{t} \int_{E_0} ((\mathbf{H}u)^2 - P_t^0(\mathbf{H}u)^2)(x)q(x)m(dx)$$
$$= \int_F u^2(x) V(dx) - \lim_{t \downarrow 0} \frac{1}{t} \int_{E_0} ((\mathbf{H}u)^2 - P_t^0(\mathbf{H}u)^2)(x)q(x)m(dx).$$

Since $q \in L^1(E_0, m)$ and $u$ is bounded, we have

(2.21) $$\frac{1}{t} \int_{E_0} ((\mathbf{H}u^2) - P_t^0(\mathbf{H}u)^2)(x)q(x)m(dx)$$
$$= \frac{1}{t} \int_{E_0} (\mathbf{H}u)^2(x)(q(x) - P_t^0 q(x))m(dx).$$

Since $q(x) = 1 - \mathbf{H}\mathbb{1}(x) = \mathbf{P}_x(\tau_0 \geq \zeta)$,

$$\begin{aligned} q(x) - P_t^0 q(x) &= \mathbf{P}_x(\tau_0 \geq \zeta) - \mathbf{P}_x(\tau_0 \geq \zeta, t < \tau_0) \\ &= \mathbf{P}_x(t \geq \tau_0 \geq \zeta) \\ &= \mathbf{E}_x \left[ \sum_{s \leq t \wedge \tau_0} \mathbb{1}_{E_0}(X_{s-}) \mathbb{1}_{\{\partial\}}(X_s) \right] \\ &= \mathbf{E}_x \left[ \int_0^{t \wedge \tau_0} \mathbb{1}_{E_0}(X_s) N(X_s, \partial) \, dH_s \right]. \end{aligned}$$

Note that $A_t := \int_0^t \mathbb{1}_{E_0}(X_s) N(X_s, \partial) \, dH_s$ is the PCAF of $X$ with Revuz measure $\mathbb{1}_{E_0}(x)\kappa(dx)$, while $t \mapsto A_{t \wedge \tau_0}$ is the PCAF of $X^0$ with the same Revuz measure $\mathbb{1}_{E_0}(x)\kappa(dx)$. Thus,

$$\begin{aligned} &\frac{1}{t} \int_{E_0} ((\mathbf{H}u^2) - P_t^0 (\mathbf{H}u)^2)(x) q(x) m(dx) \\ (2.22) \quad &= \frac{1}{t} \int_{E_0} (\mathbf{H}u)^2(x) \mathbf{E}_x[A_{t \wedge \tau_0}] m(dx) \\ &= \int_{E_0} \left( \frac{1}{t} \int_0^t P_s^0 ((\mathbf{H}u)^2)(x) \, ds \right) \kappa(dx), \end{aligned}$$

where the second equality is due to Theorem 5.1.3(iii) of [20]. On the other hand,

$$P_s^0((\mathbf{H}u)^2)(x) \leq P_s((\mathbf{H}u)^2)(x) \qquad \text{for every } s > 0 \text{ and } x \in E_0$$

and by Theorem 5.1.3(iii) and (vi) of [20],

$$\begin{aligned} &\lim_{t \downarrow 0} \int_E \left( \frac{1}{t} \int_0^t P_s((\mathbf{H}u)^2)(x) \, ds \right) \kappa(dx) \\ &= \lim_{t \downarrow 0} \frac{1}{t} \int_E (\mathbf{H}u)^2(x) \mathbf{E}_x[A_t] m(dx) \\ &= \int_{E_0} (\mathbf{H}u)^2(x) \kappa(dx) < \infty, \end{aligned}$$

since $\mathbf{H}u \in \mathcal{F}_{e,b}$. Therefore, we have, by (2.22) and the dominated convergence theorem,

$$\begin{aligned} &\lim_{t \downarrow 0} \frac{1}{t} \int_{E_0} ((\mathbf{H}u^2) - P_t^0(\mathbf{H}u)^2)(x) q(x) m(dx) \\ &= \int_{E_0} \left( \lim_{t \downarrow 0} \frac{1}{t} \int_0^t P_s^0((\mathbf{H}u)^2)(x) \, ds \right) \kappa(dx) \\ &= \int_{E_0} (\mathbf{H}u)^2(x) \kappa(dx). \end{aligned}$$

Hence, by (2.20),

$$(2.23) \qquad L(w,q) = \int_F u(x)^2 V(dx) - \int_{E_0} (\mathbf{H}u)^2(x)\kappa(dx).$$

This combined with (2.19) proves identity (2.17). □

For $\alpha > 0$, define the $\alpha$-order Feller measure $U_\alpha$ on $F \times F$ by

$$U_\alpha(f \otimes g) := \alpha(\mathbf{H}^\alpha f, \mathbf{H}g)_{L^2(E_0,m)} \qquad \text{for } f,g \in \mathcal{B}(F)_b^+.$$

Since $\mathbf{H}^\alpha f = \mathbf{H}f - \alpha G_\alpha^0 \mathbf{H}f$, it follows that

$$(2.24) \qquad \lim_{\alpha \to \infty} U_\alpha(f \otimes g) = U(f \otimes g) \qquad \text{for } f,g \in \mathcal{B}(F)_b^+.$$

We will now combine Lemma 2.4 with the next identity involving the $\alpha$-order Feller measure which first appeared in [16] in a conservative case. This identification will be utilized again in the latter part of the next section.

LEMMA 2.5. *For $\alpha > 0$ and for any bounded measurable function $u$ on $F$, let $w = \mathbf{H}(u^2) - (\mathbf{H}u)^2$. Then*

$$\alpha(\mathbf{H}^\alpha \mathbb{1}, w)_{L^2(E_0,m)} + \alpha \int_{E_0 \times F} (\mathbf{H}u(x) - u(\xi))^2 \mathbf{H}^\alpha(x, d\xi) m(dx),$$

$$= \int_{F \times F} (u(\xi) - u(\eta))^2 U_\alpha(d\xi, d\eta) + \alpha(q, \mathbf{H}^\alpha(u^2))_{L^2(E_0,m)}.$$

PROOF. Take any open set $D \subset E_0$ with $m(D) < \infty$ and put $U_D^\alpha(u,v) = \alpha(\mathbf{H}^\alpha u, \mathbf{H}v)_{L^2(D,m)}$. Then

$$\alpha(\mathbf{H}^\alpha \mathbb{1}, \mathbf{H}u^2 - (\mathbf{H}u)^2)_{L^2(D,m)} + \alpha \int_{D \times F} (\mathbf{H}u(x) - u(\xi))^2 \mathbf{H}^\alpha(x, d\xi) m(dx)$$

$$= U_D^\alpha(1, u^2) - 2U_D^\alpha(u,u) + \alpha(\mathbf{H}^\alpha u^2, 1)_{L^2(D,m)}$$

$$= \int_{F \times F} (u(\xi) - u(\eta))^2 U_D^\alpha(d\xi, d\eta) + \alpha(\mathbf{H}^\alpha u^2, q)_{L^2(D,m)}.$$

It then suffices to let $D \uparrow E_0$. □

The following theorem is an improvement of Theorem 6.1 in [19] in several respects. First, we do not assume that $X$ is conservative. Second, it gives an identity rather than an inequality, which is important in Theorem 2.7 and Corollary 2.9 below. Third, the process $X$ is allowed to have killings inside $E$.

THEOREM 2.6. *Assume condition* (2.16) *holds. For any* $u \in \mathcal{F}_e$,

$$\mu_{\langle \mathbf{H}u \rangle}(E_0) + 2 \int_{E_0 \times F} (\mathbf{H}u(x) - u(\xi))^2 J(dx, d\xi) + \int_{E_0} (\mathbf{H}u)^2(x) \kappa(dx)$$
$$= \int_{F \times F} (u(\xi) - u(\eta))^2 U(d\xi, d\eta) + 2 \int_F u(\xi)^2 V(d\xi).$$

PROOF. Without loss of generality, we may assume that $u \in \mathcal{F}_{e,b}$ since otherwise we consider $u_n = ((-n) \vee u) \wedge n$ and then pass $n \to \infty$. For $\alpha > 0$, by Lemma 2.5,

$$\text{(2.25)} \quad \alpha \int_{F \times F} (u(\xi) - u(\eta))^2 U_\alpha(d\xi, d\eta) + \alpha(q, \mathbf{H}^\alpha(u^2))_{L^2(E_0, m)}$$
$$= \alpha(\mathbf{H}^\alpha \mathbb{1}, w)_{L^2(E_0, m)} + \alpha \int_{E_0 \times F} (\mathbf{H}u(x) - u(\xi))^2 \mathbf{H}^\alpha(x, d\xi)$$

where $w = \mathbf{H}(u^2) - (\mathbf{H}u)^2$ and $q = 1 - \mathbf{H}\mathbb{1}$. It follows from (2.24) that

$$\text{(2.26)} \quad \lim_{\alpha \to \infty} \alpha \int_{F \times F} (u(\xi) - u(\eta))^2 U_\alpha(d\xi, d\eta)$$
$$= \int_{F \times F} (u(\xi) - u(\eta))^2 U(d\xi, d\eta)$$

and

$$\text{(2.27)} \quad \lim_{\alpha \to \infty} \alpha(q, \mathbf{H}^\alpha(u^2))_{L^2(E_0, m)} = \int_F u(\xi)^2 V(d\xi) \qquad \text{for every } \alpha > 0.$$

The last term in (2.25) can be written as

$$\alpha \mathbf{E}_m[e^{-\alpha \tau_0} (\mathbf{H}u(X_0) - u(X_{\tau_0}))^2 \mathbb{1}_{\{\tau_0 < \zeta\}}].$$

Since, with the exception of a set of zero capacity, every point of $F$ is a regular point of $X$ (cf. [20]) and since $u \in \mathcal{F}_{e,b}$,

$$K_t := \mathbf{H}u(X_{t \wedge \tau_0}) - \mathbf{H}u(X_0)$$

is a $\mathbf{P}_x$-martingale for q.e. $x \in E$. We claim that

$$\text{(2.28)} \quad \limsup_{t \to 0} \frac{1}{t} \mathbf{E}_m[K_t^2] < \infty.$$

This is because, by Fukushima's decomposition,

$$\mathbf{H}u(X_t) - \mathbf{H}u(X_0) = M_t^{\mathbf{H}u} + N_t^{\mathbf{H}u},$$

where $M^{\mathbf{H}u}$ is a martingale additive functional of $X$ having finite energy and $N^{\mathbf{H}u}$ is a continuous additive functional of $X$ having zero energy. Thus,

$$t \mapsto N_{t \wedge \tau_0}^{\mathbf{H}u} = K_t - M_{t \wedge \tau_0}^{\mathbf{H}u}$$

is a $\mathbf{P}_x$-martingale for q.e. $x \in E$. Since $N^{\mathbf{H}u}$ has finite energy, we have, for each fixed $t > 0$,

$$\mathbf{E}_{\mathbb{1}_{E_0}m}[[N^{\mathbf{H}u}]_t; t < \tau_0] = \mathbf{E}_{\mathbb{1}_{E_0}m}\left[\lim_{n\to\infty} \sum_{k=1}^n (N^{\mathbf{H}u}_{kt/n} - N^{\mathbf{H}u}_{(k-1)t/n})^2; t < \tau_0\right]$$

$$\leq \lim_{n\to\infty} \mathbf{E}_m\left[\sum_{k=1}^n (N^{\mathbf{H}u}_{kt/n} - N^{\mathbf{H}u}_{(k-1)t/n})^2\right]$$

$$= 0.$$

Hence, $\mathbf{P}_{\mathbb{1}_{E_0}m}$-a.s., for every $t > 0$, $[N^{\mathbf{H}u}]_t = 0$ on $\{t < \tau_0\}$. By the continuity of $[N^{\mathbf{H}u}]$, we have $[N^{\mathbf{H}u}]_{\tau_0} = 0$ $\mathbf{P}_{\mathbb{1}_{E_0}m}$-a.s. Thus, $\mathbf{P}_{\mathbb{1}_{E_0}m}$-a.s., $N^{\mathbf{H}u}_t = 0$ for every $t \leq \tau_0$. This implies that $\mathbf{P}_{\mathbb{1}_{E_0}m}$-a.s., $K_t = M^{\mathbf{H}u}_{t\wedge\tau_0}$ for every $t \leq \tau_0$. In particular,

$$\limsup_{t\to 0} \frac{1}{t}\mathbf{E}_m[K_t^2] = \limsup_{t\to 0} \frac{1}{t}\mathbf{E}_m[(M^{\mathbf{H}u}_{t\wedge\tau_0})^2] \leq [(M^{\mathbf{H}u}_t)^2] < \infty.$$

This proves the claim (2.28). However, $K$ is not a MAF of $X$ since $K_t = K_{\tau_0}$ for $t \geq \tau_0$. It is not a MAF of $X^0$ either, since

$$K_{\tau_0} - K_{\tau_0-} = \mathbf{H}u(X_{\tau_0}) - \mathbf{H}u(X_{\tau_0-}) = (u(X_{\tau_0}) - \mathbf{H}u(X_{\tau_0-}))\mathbb{1}_{\{X_{\tau_0}\neq X_{\tau_0-}\}}$$

is not a function of $X_{\tau_0-}$ on $\{X_{\tau_0} \neq X_{\tau_0-}\}$. However, if we define

$$M_t = \begin{cases} \mathbf{H}u(X_t) - \mathbf{H}u(X_0), & \text{for } 0 \leq t < \tau_0, \\ \mathbf{H}u(X_{\tau_0-}) - \mathbf{H}u(X_0), & \text{for } t \geq \tau_0, \end{cases}$$

then $M$ is a MAF of $X^0$ having finite energy. We will use $\langle M \rangle$ to denote the dual predictable projection of the variational process $[M]$ of $M$. Clearly, $\langle M \rangle$ is a PCAF of $X^0$ and we will use $\mu_{\langle M \rangle}$ to denote the Revuz measure of $\langle M \rangle$ with respect to the subprocess $X^0$. Note that $\tau_0$ is the lifetime of $X^0$. By Itô's formula,

$$dM_t^2 = 2M_{t-}\, dM_t + d[M]_t$$

and so

$$d(e^{-\alpha t} M_t^2) = -\alpha e^{-\alpha t} M_{t-}^2\, dt + 2\alpha e^{-\alpha t} M_{t-}\, dM_t + e^{-\alpha t}\, d[M]_t.$$

Using dual predictable projection, we have

$$\alpha \mathbf{E}_m[e^{-\alpha\tau_0} M_{\tau_0}^2]$$

$$= \alpha \mathbf{E}_{\mathbb{1}_{E_0}m}\left[-\alpha \int_0^\infty e^{-\alpha t} M_{t-}^2 \mathbb{1}_{\{t<\tau_0\}}\, dt + \int_0^{\tau_0} e^{-\alpha t}\, d[M]_t\right]$$

$$= \alpha \mathbf{E}_{\mathbb{1}_{E_0}m}\left[-\alpha \int_0^\infty e^{-\alpha t} \langle M \rangle_t \mathbb{1}_{\{t<\tau_0\}}\, dt + \int_0^\infty e^{-\alpha t}\, d\langle M \rangle_t\right]$$

$$= \alpha \mathbf{E}_{1_{E_0} m} \left[ -\int_0^\infty e^{-s} \langle M \rangle_{s/\alpha} 1_{\{s/\alpha < \tau_0\}} \, ds + \int_0^\infty e^{-\alpha t} \, d\langle M \rangle_t \right]$$

$$= \alpha \mathbf{E}_{1_{E_0} m} \left[ -\int_0^\infty e^{-s} \langle M \rangle_{s/\alpha} \, ds + \int_0^\infty e^{-\alpha t} \, d\langle M \rangle_t \right]$$

$$+ \alpha \mathbf{E}_{1_{E_0} m} \left[ \int_0^\infty e^{-s} \langle M \rangle_{s/\alpha} 1_{\{s/\alpha \geq \tau_0\}} \, ds \right].$$

Consequently,

$$\lim_{\alpha \to \infty} \alpha \mathbf{E}_m [e^{-\alpha \tau_0} M_{\tau_0}^2]$$

$$= -\int_0^\infty s e^{-s} \mu_{\langle M \rangle}(E_0) \, ds + \mu_{\langle M \rangle}(E_0)$$

(2.29)
$$+ \lim_{\alpha \to \infty} \alpha \mathbf{E}_{1_{E_0} m} \left[ \int_0^\infty e^{-s} \langle M \rangle_{s/\alpha} 1_{\{s/\alpha \geq \tau_0\}} \, ds \right]$$

$$= \int_0^\infty e^{-s} \left( \lim_{\alpha \to \infty} \alpha \mathbf{E}_{1_{E_0} m} [\langle M \rangle_{s/\alpha} 1_{\{s/\alpha \geq \tau_0\}}] \right) ds$$

$$= 0,$$

where in the second to the last equality we used the dominated convergence theorem and in the last equality we applied Lemma 2.1 with $X^0$ in place of $X$. Next, note that, by using the Lévy system of $X$, we have

$$\mathbf{E}_m [e^{-\alpha \tau_0} (\mathbf{H} u(X_{\tau_0 -}) - u(X_{\tau_0}) 1_{\{\tau_0 < \zeta\}})^2]$$

$$= \mathbf{E}_m [e^{-\alpha \tau_0} (\mathbf{H} u(X_{\tau_0 -}) - u(X_{\tau_0}) 1_{\{X_{\tau_0 -} \neq X_{\tau_0}\}})^2 1_{\{\tau_0 < \zeta\}}]$$

$$= \mathbf{E}_{1_{E_0} m} \left[ \int_0^{\tau_0} e^{-\alpha s} \int_F (\mathbf{H} u(x) - u(\xi))^2 N(X_s, d\xi) \, dH_s \right]$$

$$= \int_{E_0} G_\alpha^0 1(x) \left( \int_F (\mathbf{H} u(x) - u(\xi))^2 N(x, d\xi) \right) \mu_H(dx).$$

Hence, by the monotone convergence theorem,

$$\lim_{\alpha \to \infty} \alpha \mathbf{E}_m [e^{-\alpha \tau_0} (\mathbf{H} u(X_{\tau_0 -}) - u(X_{\tau_0}))^2 1_{\{\tau_0 < \zeta\}}]$$

$$= \int_{E_0} \left( \int_F (\mathbf{H} u(x) - u(\xi))^2 N(x, d\xi) \right) \mu_H(dx)$$

$$= \int_{E_0 \times F} (\mathbf{H} u(x) - u(\xi))^2 N(x, d\xi) \mu_H(dx)$$

$$= 2 \int_{E_0 \times F} (\mathbf{H} u(x) - u(\xi))^2 J(dx, d\xi).$$

It then follows from (2.29) that

$$\lim_{\alpha \to \infty} \alpha \mathbf{E}_m [e^{-\alpha \tau_0} (\mathbf{H} u(X_0) - u(X_{\tau_0}))^2 1_{\{\tau_0 < \zeta\}}]$$

$$= \lim_{\alpha\to\infty} \alpha\mathbf{E}_m[e^{-\alpha\tau_0}(-M_{\tau_0} + \mathbf{H}u(X_{\tau_0-}) - u(X_{\tau_0}))^2 \mathbb{1}_{\{\tau_0<\zeta\}}]$$

$$= \lim_{\alpha\to\infty} \alpha\mathbf{E}_m[e^{-\alpha\tau_0}(\mathbf{H}u(X_{\tau_0-}) - u(X_{\tau_0}))^2 \mathbb{1}_{\{\tau_0<\zeta\}}]$$

$$= 2\int_{E_0\times F}(\mathbf{H}u(x) - u(\xi))^2 J(dx,d\xi).$$

Passing $\alpha\to\infty$ in identity (2.25), the above calculation together with (2.17) and (2.26) and (2.27) proves the theorem. □

THEOREM 2.7. *Assume condition (2.16) holds. For any $u\in\mathcal{F}_e$,*

$$\mathcal{E}(\mathbf{H}u,\mathbf{H}u) = \tfrac{1}{2}\mu^c_{\langle\mathbf{H}u\rangle}(F) + \int_{F\times F}(u(x)-u(y))^2(\tfrac{1}{2}U(dx,dy) + J(dx,dy))$$
$$+ \int_F u(x)^2(V(dx)+\kappa(dx)).$$

PROOF. It is well known (cf. [20]) that

$$\mathcal{E}(\mathbf{H}u,\mathbf{H}u)$$
$$= \tfrac{1}{2}\mu_{\langle\mathbf{H}u\rangle}(E) + \tfrac{1}{2}\mu^k_{\langle\mathbf{H}u\rangle}(E)$$
$$= \tfrac{1}{2}\mu_{\langle\mathbf{H}u\rangle}(E_0) + \tfrac{1}{2}\mu_{\langle\mathbf{H}u\rangle}(F) + \tfrac{1}{2}\mu^k_{\langle\mathbf{H}u\rangle}(E)$$
$$= \tfrac{1}{2}\mu_{\langle\mathbf{H}u\rangle}(E_0) + \tfrac{1}{2}\mu^k_{\langle\mathbf{H}u\rangle}(E_0) + \tfrac{1}{2}\mu^c_{\langle\mathbf{H}u\rangle}(F) + \tfrac{1}{2}\mu^j_{\langle\mathbf{H}u\rangle}(F) + \mu^k_{\langle\mathbf{H}u\rangle}(F),$$

where

$$\mu^j_{\langle\mathbf{H}u\rangle}(dx) = \left(\int_E(\mathbf{H}u(x)-\mathbf{H}u(y))^2 N(x,dy)\right)\mu_H(dx)$$

and

$$\mu^k_{\langle\mathbf{H}u\rangle}(dx) = (\mathbf{H}u)^2(x)\kappa(dx) = u(x)^2\kappa(dx).$$

The assertion of the theorem now follows from Theorem 2.6. □

We now study the trace of the process $X$ on the quasi-closed set $F$. A quasi-support of a Borel measure is the smallest quasi-closed set outside which the measure has zero charge. The quasi-support is unique up to quasi-equivalence. Denote by $\overset{\circ}{S}$ the family of all positive Radon measures on $E$ charging no set of zero capacity. Put

(2.30) $\qquad \overset{\circ}{S}_F = \{\mu\in\overset{\circ}{S}:$ the quasi support of $\mu = F$ q.e.$\}.$

LEMMA 2.8. *For a quasi-closed subset $F\subset E$ with positive capacity, $\overset{\circ}{S}_F\neq\varnothing$.*

PROOF. The proof is the same as that for Lemma 4.1 in [19]. But for the reader's convenience, we spell out the details here. Let $\sigma_F := \inf\{t \geq 0 : X_t \in F\}$ and $g$ be a strictly positive function in $L^1(E, m)$. Define, for $B \in \mathcal{B}(E)$,

$$\mu(B) := \int_E g(x) \mathbf{P}_x(X_{\sigma_F} \in B, \sigma_F < \infty) m(dx).$$

Since $F$ is quasi-closed and has positive capacity, $\mu(E \setminus F) = 0$ and $\mu$ is a nontrivial finite measure charging no set of zero capacity. If $f \in \mathcal{F}$ and $f = 0$ $\mu$-a.e on $E$, then

$$\int_E g(x) \mathbf{H}|f|(x) m(dx) = \int_E |f(x)| \mu(dx) = 0.$$

Since $\mathbf{H}|f|$ is quasi-continuous, $\mathbf{H}|f| = 0$ q.e. on $E$ and, thus, $f = 0$ q.e. on $F$. It now follows from [20], Theorem 4.6.2, that $F$ is a quasi-support of $\mu$; that is, $\mu \in \overset{\circ}{S}_F$. The lemma is proved. $\square$

Fix a measure $\mu \in \overset{\circ}{S}_F$, and let $A^\mu$ be the PCAF of $X$ with Revuz measure $\mu$. Let $N$ be the exceptional set in the definition of $A^\mu$. Then the support of $A^\mu$, which is defined as (see (5.1.21) of [20])

$$\{x \in E \setminus N : \mathbf{P}_x(\inf\{t > 0 : A_t^\mu > 0\} = 0) = 1\},$$

is nearly Borel, finely closed and equals to $F$ q.e. Therefore, we may and shall assume that the support of $A^\mu$ is just $F$ in accordance with the preceding convention. Note that $\mathbf{H}u$ with $u \in \mathcal{F}_e$ remains the same q.e. if $F$ is replaced by another quasi-closed set that equals to $F$ q.e.

Let $Y$ be the time-changed process of the Hunt process $X$ by the right continuous inverse of $A^\mu$. Then by [30], (65.9), $Y$ is a $\mu$-symmetric right process on $F$. The Dirichlet form on $(\check{\mathcal{E}}, \check{\mathcal{F}})$ of $Y$ on $L^2(F; \mu)$ admits an expression (cf. Theorem 6.2.1 of [20]):

$$\check{\mathcal{F}}_e = \mathcal{F}_e|_F,$$
(2.31) $$\check{\mathcal{F}} = \mathcal{F}_e|_F \cap L^2(F, \mu),$$
$$\check{\mathcal{E}}(u, v) = \mathcal{E}(\mathbf{H}u, \mathbf{H}v) \qquad \text{for } u, v \in \check{\mathcal{F}}_e.$$

By [14] and [28], the Dirichlet form $(\check{\mathcal{E}}, \check{\mathcal{F}})$ is quasi-regular on $F$. See the Appendix of this paper for the existence of a Lévy system for $Y$ and the probabilistic characterization for the Beurling–Deny decomposition of $(\check{\mathcal{E}}, \check{\mathcal{F}})$. Recall our convention that functions in $\mathcal{F}_e$ are always represented by their quasi-continuous versions. It can also be verified that the extended Dirichlet space $\check{\mathcal{F}}_e$ of $(\check{\mathcal{F}}, \check{\mathcal{E}})$ coincides with $\mathcal{F}_e|_F$, independent of the choice of $\mu \in \overset{\circ}{S}_F$.

Denote by $\widetilde{F}$ the topological support of the measure $\mu$. The trace Dirichlet form $(\check{\mathcal{E}}, \check{\mathcal{F}})$, which is quasi-regular, is, in general, not a regular Dirichlet form on $F$. However, by [20], Theorem 6.2.1, it is regular on $\widetilde{F}$. Here we note that $\mathrm{Cap}(F \setminus \widetilde{F}) = 0$, $\mu(\widetilde{F} \setminus F) = 0$ and so $L^2(F;\mu) = L^2(\widetilde{F};\mu)$. But it is possible that $\mathrm{Cap}(\widetilde{F} \setminus F) > 0$ (see [20], Example 5.1.2). However, $\widetilde{F} \setminus F$ always has zero $\check{\mathcal{E}}_1$-capacity by [20], Theorem 6.2.1(iv). Hence, the natural inclusion map $i: F \to \widetilde{F}$ gives a quasi-homeomorphism between the quasi-regular Dirichlet form $(\check{\mathcal{E}}, \check{\mathcal{F}})$ on $L^2(F;\mu)$ and the regular Dirichlet form $(\check{\mathcal{E}}, \check{\mathcal{F}})$ on $L^2(\widetilde{F};\mu)$. We shall employ this useful fact in the proof of the next theorem. There are some more related discussions on this at the end of the next section.

By Theorem 2.7, we have the following:

COROLLARY 2.9. *Assume condition* (2.16) *holds. For* $u \in \check{\mathcal{F}}_e$,

$$\check{\mathcal{E}}(u,u) = \tfrac{1}{2}\mu^c_{\langle \mathbf{H}u \rangle}(F) + \int_{F \times F} (u(x) - u(y))^2 (\tfrac{1}{2} U(dx,dy) + J(dx,dy))$$
$$+ \int_F u(x)^2 (V(dx) + \kappa(dx)). \tag{2.32}$$

We now show that the decomposition in (2.32) is the Beurling–Deny decomposition for $(\check{\mathcal{E}}, \check{\mathcal{F}})$.

THEOREM 2.10. *Assume condition* (2.16) *holds. The bilinear form* $(u,v) \mapsto \mu^c_{\langle \mathbf{H}u, \mathbf{H}v \rangle}(F)$ *has the strong local property on* $\check{\mathcal{F}}$; *that is, if* $u, v \in \check{\mathcal{F}} \cap C_c(\widetilde{F})$ *and* $u$ *is constant in a neighborhood of* $\mathrm{supp}[v]$, *then* $\mu^c_{\langle \mathbf{H}u, \mathbf{H}v \rangle}(F) = 0$. *In other words,*

$$\check{\mathcal{E}}^c(u,v) = \tfrac{1}{2}\mu^c_{\langle \mathbf{H}u, \mathbf{H}v \rangle}(F) \qquad \text{for } u, v \in \check{\mathcal{F}}.$$

PROOF. As is mentioned above, $(\check{\mathcal{E}}, \check{\mathcal{F}})$ is a regular Dirichlet form on $L^2(\widetilde{F}, \mu)$. Let $u \in \check{\mathcal{F}}_b \cap C_c(\widetilde{F})$ such that $u = c$ for some $c \in \mathbf{R}$ on a relative open subset $I$ of $\widetilde{F}$. For any relatively compact open subset $I_1$ of $I$, there are open subsets $D_1$ and $D$ of $E$ such that $\overline{D}_1 \subset D$, $\overline{D}$ is compact, $D_1 \cap F = I_1$ and $D \cap F = I$. Since $(\mathcal{E}, \mathcal{F})$ is a regular Dirichlet space on $L^2(E,m)$, there is a $\varphi \in \mathcal{F} \cap C_c(E)$ such that $\varphi = 1$ on $D_1$ and $\varphi = 0$ on $D^c$. Let $v = c\varphi + (1 - \varphi)\mathbf{H}u$. Then $v \in \mathcal{F}_{e,b}$ and $v$ is constant on $D_1$. Hence, $\mu^c_{\langle v \rangle}(D_1) = 0$ and, thus, we conclude

$$\mu^c_{\langle v \rangle}(I_1) = 0. \tag{2.33}$$

Since $v = u$ q.e. on $F$, we have $\mathbf{H}v = \mathbf{H}u$ q.e. Define $v_0 = v - \mathbf{H}v$, which is in $\mathcal{F}^0_{e,b}$. Let $\{K_n, n \geq 1\}$ be a generalized nest so that $F \cap K_n$ is a closed set. Let $(\mathcal{E}, \mathcal{F}^{E \setminus (F \cap K_n)})$ be the Dirichlet space for the subprocess of $X$ killed upon leaving $E \setminus (F \cap K_n)$. Clearly, $v_0 \in \mathcal{F}^0_{e,b} \subset \mathcal{F}^{E \setminus (F \cap K_n)}_{e,b}$. Since $(\mathcal{E}, \mathcal{F}^{E \setminus (F \cap K_n)})$ is regular on $L^2(E \setminus (F \cap K_n); m)$, by Theorem 4.4.3 of [20], $\mu^c_{\langle \psi \rangle}(F \cap K_n) = 0$ for $\psi \in \mathcal{F} \cap C_c(E \setminus (F \cap K_n))$ and, hence, for $\psi = v_0$. Thus, $\mu^c_{\langle \psi \rangle}(F \cap K_n) = 0$ and passing $n$ to infinity we have $\mu^c_{\langle \psi \rangle}(F) = 0$. In particular, $\mu^c_{\langle v_0 \rangle}(I \cap F) = 0$. On the other hand,

$$\begin{aligned}\mu^c_{\langle v \rangle}(I_1 \cap F) &= \mu^c_{\langle \mathbf{H}v + v_0 \rangle}(I_1 \cap F) \\ &= \mu^c_{\langle \mathbf{H}v \rangle}(I_1 \cap F) + 2\mu^c_{\langle \mathbf{H}v, v_0 \rangle}(I \cap F) + \mu^c_{\langle v_0 \rangle}(I_1 \cap F) \\ &= \mu^c_{\langle \mathbf{H}v \rangle}(I_1 \cap F) \\ &= \mu^c_{\langle \mathbf{H}u \rangle}(I_1 \cap F).\end{aligned}$$

Thus, by (2.33), $\mu^c_{\langle \mathbf{H}u \rangle}(I_1 \cap F) = 0$. By letting $I_1 \uparrow I$, we get

(2.34) $$\mu^c_{\langle \mathbf{H}u \rangle}(I \cap F) = 0.$$

Now for $u, v \in \check{\mathcal{F}} \cap C_c(\widetilde{F})$ such that $u$ is constant in a neighborhood of $\operatorname{supp}[v]$, we let $F_1 = \operatorname{supp}[v]$ and $F_2 = \widetilde{F} \setminus \operatorname{supp}[v]$. By (2.34),

$$\mu^c_{\langle \mathbf{H}u \rangle}(F_1 \cap F) = 0 \quad \text{and} \quad \mu^c_{\langle \mathbf{H}v \rangle}(F_2 \cap F) = 0.$$

Since $F \subset \widetilde{F}$ q.e. and $\mu^c_{\langle \mathbf{H}u, \mathbf{H}v \rangle}$ does not charge on sets of zero capacity, it follows then

$$\begin{aligned}|\mu^c_{\langle \mathbf{H}u, \mathbf{H}v \rangle}(F)| &= |\mu^c_{\langle \mathbf{H}u, \mathbf{H}v \rangle}(F_1 \cap F) + \mu^c_{\langle \mathbf{H}u, \mathbf{H}v \rangle}(F_2 \cap F)| \\ &\leq \sqrt{\mu^c_{\langle \mathbf{H}u \rangle}(F_1 \cap F) \mu^c_{\langle \mathbf{H}v \rangle}(F_1 \cap F)} \\ &\quad + \sqrt{\mu^c_{\langle \mathbf{H}u \rangle}(F_2 \cap F) \mu^c_{\langle \mathbf{H}v \rangle}(F_2 \cap F)} \\ &= 0.\end{aligned}$$

This proves the theorem. □

It follows from Corollary 2.9 and Theorem 2.10 that

(2.35) $$\check{J}(dx, dy) := \tfrac{1}{2} U(dx, dy) + J(dx, dy)|_{F \times F}$$

and

(2.36) $$\check{\kappa}(dx) := V(dx) + \kappa(dx)|_F$$

are the jumping measure and the killing measure for the time-changed process $Y$, respectively.

We now proceed to remove condition (2.16).

THEOREM 2.11. *Theorem 2.7, Corollary 2.9 and Theorem 2.10 remain true without the assumption* (2.16).

PROOF. The only place where condition (2.16) is used is to ensure that the symmetry can be applied in (2.18) and in (2.21). By using a time-change method, the assumption (2.16) can be dropped. Let $\varphi > 0$ on $E_0$ be such that $\int_{E_0} \varphi(x) m(dx) < \infty$. Let

$$\nu(dx) := (\mathbb{1}_{E_0}(x)\varphi(x) + \mathbb{1}_{E_0^c}(x))m(dx).$$

Let $A^\nu$ be the PCAF of $X$ with Revuz measure $\nu$ and

$$\sigma_t := \inf\{s \geq 0 : A_s^\nu > t\}.$$

Then $Z_t := X_{\sigma_t}$ is a $\nu$-symmetric right process on $E$ with Dirichlet form $(\mathcal{E}^Z, \mathcal{F}^Z)$. It is well known (cf. [20]) that $\mathcal{F}_e^Z = \mathcal{F}_e$, $\mathcal{E}^Z = \mathcal{E}$ on $\mathcal{F}_e$. Clearly, $\mathbf{H}^Z u = \mathbf{H} u$ for $u \in \mathcal{F}_e$. Since $\nu(E_0) < \infty$, Theorem 2.7, Corollary 2.9 and Theorem 2.10 hold for the process $Z$. Hence, Theorem 2.7 and Corollary 2.9 hold without the assumption (2.16). Moreover, (2.32) (for $Z$) gives the Beurling–Deny decomposition for the Dirichlet form of the time changed process of $X$ by $\mu \in \mathring{S}_F$.

We now show that Feller measures are invariant under time changes. For this, let $U^X$, $U^Z$ and $L^X$, $L^Z$ denote the Feller measures and energy functionals for the processes $X$ and $Z$, respectively. We will use $G^{X,0}$ and $G^{Z,0}$ to denote the 0-resolvent of the part process $X^0$ and $Z^0$ in $E_0$, respectively. Clearly, for $f \geq 0$ on $E_0$ and $x \in E_0$,

$$G^{Z,0} f(x) = \mathbf{E}_x\left[\int_0^\infty f(X_{\sigma_t}^0)\, dt\right]$$

$$= \mathbf{E}_x\left[\int_0^\infty f(X_t^0)\varphi(X_t)\, dt\right]$$

$$= G^{X,0}(\varphi f)(x).$$

By (2.14) above and Proposition 3.6(ii) in [21], for $u, v \in \mathcal{B}(F)_b^+$,

$$U^X(u \otimes v) = L^X(\mathbf{H}u, \mathbf{H}v)$$

$$= \sup\left\{\int_{E_0} \mathbf{H}u(x)\mu(dx) : (\mu G^{X,0})(dx) \leq \mathbf{H}v(x)m(dx) \text{ on } E_0\right\}.$$

Hence,

$$U^Z(u \otimes v)$$

$$= \sup\left\{\int_{E_0} \mathbf{H}u(x)\mu(dx) : (\mu G^{Z,0})(dx) \leq \mathbf{H}v(x)\nu(dx) \text{ on } E_0\right\}$$

$$= \sup\left\{\int_{E_0} \mathbf{H}u(x)\mu(dx) : \int_{E_0} G^{Z,0}f(x)\mu(dx) \leq \int_{E_0} f(x)\mathbf{H}v(x)\nu(dx)\right.$$

$$\left. \text{for every } f \geq 0 \text{ on } E_0\right\}$$

$$\begin{aligned}
&= \sup\biggl\{\int_{E_0} \mathbf{H}u(x)\mu(dx) : \int_{E_0} G^{X,0}(\varphi f)(x)\mu(dx) \\
&\qquad\qquad\qquad\qquad \le \int_{E_0} f(x)\mathbf{H}v(x)\varphi(x)m(dx) \\
&\qquad\qquad\qquad\qquad \text{for every } f \ge 0 \text{ on } E_0\biggr\} \\
&= \sup\biggl\{\int_{E_0} \mathbf{H}u(x)\mu(dx) : \int_{E_0} G^{X,0}\mu(x)(\varphi f)(x)m(dx) \\
&\qquad\qquad\qquad\qquad \le \int_{E_0} (f\varphi)(x)\mathbf{H}v(x)m(dx) \\
&\qquad\qquad\qquad\qquad \text{for every } f \ge 0 \text{ on } E_0\biggr\} \\
&= \sup\biggl\{\int_{E_0} \mathbf{H}u(x)\mu(dx) : G^{X,0}\mu \le \mathbf{H}v \text{ on } E_0\biggr\} \\
&= U^X(u \otimes v).
\end{aligned}$$

This shows that $U^X = U^Z$ and, consequently, $V^X = V^Z$. In particular, this shows that $U^X$ is a symmetric measure on $F \times F$. Therefore, Theorem 2.10 holds for the process $X$ without the assumption (2.16). $\square$

EXAMPLE 2.12. Let $X = (X_t, \mathbf{P}_x)_{x \in \mathbb{R}^n}$ be the standard Brownian motion on $\mathbb{R}^n$ with $n \ge 3$. Let $S$ be a $C^3$ compact hypersurface so that $E_0 = \mathbb{R}^n \setminus S$ is the union of the interior domain $D_i$ and exterior domain $D_e$. We denote by $\sigma$ the surface measure on $S$. Further, $\frac{\partial}{\partial n^i_\xi}$ and $\frac{\partial}{\partial n^e_\xi}$ will denote the inward normal and outward normal derivative at $\xi \in S$ from the view of $D_i$, respectively. We consider the Poisson kernel $K(x, \xi), x \in E_0, \xi \in S$, and the escape probability of $X$ from $S$ defined by

$$q(x) = 1 - \mathbf{H}\mathbb{1}(x) = \mathbf{P}_x(\sigma_S = \infty), \qquad x \in E_0.$$

Note that $q(x) > 0$ only for $x \in D_e$.

Then the Feller measure $U$ and the supplementary Feller measure $V$ with respect to $X$ and $S$ have densities $U(\xi, \eta)$ and $v(\xi)$ relative to $\sigma \times \sigma$ and $\sigma$, respectively, which admit the expressions

$$(2.37) \qquad U(\xi, \eta) = \frac{1}{2}\frac{\partial K(\xi, \eta)}{\partial n^i_\xi} + \frac{1}{2}\frac{\partial K(\xi, \eta)}{\partial n^e_\xi}, \qquad \xi \ne \eta, \xi, \eta \in S,$$

$$(2.38) \qquad v(\xi) = \frac{1}{2}\frac{\partial q(\xi)}{\partial n^e_\xi}, \qquad \xi \in S.$$

Formula (2.37) has been established in [19], Example 2.1, where it is also shown that, in the special case of $S = \Sigma_r$, the sphere of radius $r > 0$ centered at the origin, (2.37) is reduced to

$$U(\xi, \eta) = \frac{2}{\Omega_n} |\xi - \eta|^{-n}, \qquad \xi, \eta \in \Sigma_r, \xi \neq \eta,$$

where $\Omega_n$ denotes the area of the unit sphere in $\mathbb{R}^n$. In this special case, (2.38) is also reduced to a constant function

$$v(\xi) = \frac{1}{2} \lim_{|x| \downarrow r} \frac{q(x)}{|x| - r} = \frac{n-2}{2r}, \qquad \xi \in \Sigma_r,$$

because $q(x) = 1 - \frac{r^{n-2}}{|x|^{n-2}}$ for $|x| > r$.

For the proof of (2.37) and (2.38), it suffices to show the formula

$$\text{(2.39)} \quad \tfrac{1}{2} \int_{\mathbb{R}^n} |\nabla \mathbf{H} f(x)|^2 \, dx = \tfrac{1}{2} \int_{S \times S} (f(\xi) - f(\eta))^2 U(\xi, \eta) \sigma(d\xi) \sigma(d\eta) + \int_S f(\xi)^2 v(\xi) \sigma(d\xi)$$

holding for any $f \in C_c^\infty(\mathbb{R}^n)$.

Take a sufficiently large $R$ such that the ball $B_R = \{x \in \mathbb{R}^n : |x| < R\}$ contains the surface $S$. Let $X^R$ be the reflecting Brownian motion on $\overline{B}_R$, whose Dirichlet form $(\mathcal{E}^R, \mathcal{F}^R)$ on $L^2(B_R)$ is given by

$$\mathcal{E}^R(u, v) = \tfrac{1}{2} \int_{B_R} \nabla u \cdot \nabla v \, dx \quad \text{and} \quad \mathcal{F}^R = H^1(B_R).$$

Here $H^1(B_R)$ is the space of $L^2$-integrable functions on $B_R$ with $L^2$-integrable first-order distributional derivatives on $B_R$. Let $U^R$ be the Feller measure on the closed set $F = S \cup \Sigma_R (\subset \overline{B}_R)$ relative to $X^R$. The surface measure on $F$ is denoted by $\sigma$ again. Let $D_{e,R} := D_e \cap B_R$. Just as in [20], Example 2.1, we see that $U^R$ admits a density with respect to $\sigma \times \sigma$ with

$$\text{(2.40)} \quad U^R(\xi, \eta) = \frac{1}{2} \frac{\partial K(\xi, \eta)}{\partial n_\xi^i} I_S(\eta) + \frac{1}{2} \frac{\partial K^R(\xi, \eta)}{\partial n_\xi^e}, \qquad \xi \in S, \eta \in F,$$

where $K(x, \eta), x \in D_i, \eta \in S$, is, as before, the Poisson kernel in $D_i$, while $K^R(x, \eta), x \in D_{e,R}, \eta \in F$, is the Poisson kernel in the region $D_{e,R}$ so that, for $x \in D_{e,R}$,

$$\text{(2.41)} \quad \int_\Gamma K^R(x, \eta) \sigma(d\eta) = \begin{cases} \mathbf{P}_x(\sigma_S < \sigma_{\Sigma_R}, X_{\sigma_S} \in \Gamma), & \Gamma \subset S, \\ \mathbf{P}_x(\sigma_S > \sigma_{\Sigma_R}, X_{\sigma_{\Sigma_R}} \in \Gamma), & \Gamma \subset \Sigma_R, \end{cases}$$

where $X$ is Brownian motion in $\mathbb{R}^n$.

Since $X^R$ is conservative and $B_R$ is of finite Lebesgue measure, Theorem 6.2 of [19] implies that

$$\tfrac{1}{2}\int_{B_R}|\nabla \mathbf{H}f(x)|^2\,dx = \tfrac{1}{2}\int_{F\times F}(f(\xi)-f(\eta))^2 U^R(d\xi,d\eta), \qquad f\in C_c^\infty(\mathbb{R}^n).$$

Therefore, if we take for any $f\in C_c^\infty(\mathbb{R}^n)$ a large $R$ such that the support of $f$ is contained in $B_R$, then we have

$$\begin{aligned}\tfrac{1}{2}\int_{B_R}|\nabla \mathbf{H}f(x)|^2\,dx &= \tfrac{1}{2}\int_{S\times S}(f(\xi)-f(\eta))^2 U^R(\xi,\eta)\sigma(d\xi)\sigma(d\eta)\\ &\quad + \int_S f(\xi)^2 v^R(\xi)\sigma(d\xi),\end{aligned} \qquad (2.42)$$

where

$$v^R(\xi) = \int_{\Sigma_R} U^R(\xi,\eta)\sigma(d\eta).$$

In view of (2.40) and (2.41), we have

$$(2.43)\qquad v^R(\xi) = \frac{1}{2}\frac{\partial q^R(\xi)}{\partial n_\xi^e}$$

for $\xi\in S$, where $q^R(x) = \mathbf{P}_x(\sigma_{\Sigma_R} < \sigma_S)$ for $x\in D_{e,R}$.

We now prove that

$$(2.44)\qquad U^R(\xi,\eta)\uparrow U(\xi,\eta) \quad \text{and}\quad \frac{\partial q^R(\xi)}{\partial n_\xi^e}\downarrow \frac{\partial q(\xi)}{\partial n_\xi^e} \qquad \text{as } R\uparrow\infty.$$

Then the desired identity (2.39) follows from (2.42) and (2.43).

For $\xi\in S$, let $K_\alpha(x,\xi), x\in D_e$, and $K_\alpha^R(x,\xi), x\in D_{e,R}$, be the $\alpha$-order Poisson kernels on $D_e$ and $D_{e,R}$, respectively. As $R\uparrow\infty$, $K^R(x,\xi)$ [respectively, $K_\alpha^R(x,\xi)$] increases to $K(x,\xi)$ [resp. $K_\alpha(x,\xi)$]. On the other hand, we see from [19], Theorem 6.2, that, as $\alpha\uparrow\infty$,

$$\alpha\int_{D_{e,R}} K^R(x,\xi)K_\alpha^R(x,\eta)\,dx \uparrow \frac{1}{2}\frac{\partial K^R(\xi,\eta)}{\partial n_\xi^e}, \qquad \xi,\eta\in S$$

$$\alpha\int_{D_e} K(x,\xi)K_\alpha(x,\eta)\,dx \uparrow \frac{1}{2}\frac{\partial K(\xi,\eta)}{\partial n_\xi^e}, \qquad \xi,\eta\in S.$$

Hence, by interchanging the order of taking limits, we get the first part of (2.44).

Finally, we take $R_1$ with $S\subset B_{R_1}$ and denote by $\sigma_1$ the surface measure on $\Sigma_{R_1}$. Then we have

$$q^R(x) = \int_{\Sigma_{R_1}} K^{R_1}(x,\eta) q^R(\eta)\sigma_1(d\eta), \qquad R_1 < R,$$

and

$$q(x) = \int_{\Sigma_{R_1}} K^{R_1}(x,\eta) q(\eta) \sigma_1(d\eta).$$

By taking the outward normal derivative in $x$ at $S$ on both sides of each of the above two equations, we arrive at the second part of (2.44) because $q^R(\eta)$ decreases to $q(\eta)$ as $R \uparrow \infty$ for each $\eta \in \Sigma_{R_1}$.

**3. Space of functions with finite Douglas integrals.** In the preceding section we have established a Beurling–Deny decomposition of the trace Dirichlet form $(\check{\mathcal{E}}, \widetilde{\mathcal{F}})$ with the jumping and killing measures being given by (1.6). This particularly means that, for any function $\varphi$ in the (extended) trace Dirichlet space $\check{\mathcal{F}}_e$, the value of the integral (1.7) is finite. We may call this value the (generalized) *Douglas integral* of the function $\varphi$ (cf. [19]).

In this section we shall look for conditions to guarantee the coincidence of the (extended) trace Dirichlet space $\check{\mathcal{F}}_e$ with the space of functions on $F$ with finite Douglas integrals. For this purpose, we shall first study the relationship between the Dirichlet form $(\mathcal{E}, \mathcal{F})$ of the given process $X$ on $E$ and the reflected Dirichlet space $(\mathcal{F}^{\mathrm{ref}}, \mathcal{E}^{\mathrm{ref}})$ for the Dirichlet form $(\mathcal{E}^0, \mathcal{F}^0)$ of the absorbed process $X^0$ on $E_0 = E \setminus F$. The notion of the reflected Dirichlet space was introduced by Silverstein (see [31]) for a transient regular Dirichlet form and was further studied in detail by the first author [6].

We continue to work with a regular irreducible Dirichlet form $(\mathcal{E}, \mathcal{F})$ on $L^2(E, m)$, an associated Hunt process $X = \{\Omega, X_t, \zeta, P_x\}_{x \in E}$ on $E$ and a quasi-closed set $F$ of $E$ satisfying (2.1). As in Section 2, each element of the space $\mathcal{F}$ will be represented by its quasi-continuous version and we will assume, without loss of generality, that $F$ is nearly Borel and finely closed. Then $E_0 = E \setminus F$ is finely open and we have by Lemma 2.2 that the Dirichlet form $(\mathcal{E}^0, \mathcal{F}^0)$ on $L^2(E_0; m)$, defined by (2.3), is transient and quasi-regular. But $(\mathcal{E}^0, \mathcal{F}^0)$ is, in general, not regular on $L^2(E_0; m)$. The first part of this section is to generalize the notion of the reflected Dirichlet space from the regular Dirichlet form setting (see [6]) to the quasi-regular Dirichlet form $(\mathcal{E}^0, \mathcal{F}^0)$ by using the energy functional $L$ defined by (2.9) for $X^0$.

Recall that $\tau_0$ is defined by (2.2), which denotes the first exit time from $E_0 = E \setminus F$ by $X$. The process $X^0$ can then be realized as

$$X^0 = \{\Omega, \mathcal{M}^0, X_t^0, \zeta^0, \mathbf{P}_x\}_{x \in E_0},$$

where

$$\zeta^0 = \tau_0 \quad \text{and} \quad X_t^0 = \begin{cases} X_t, & \text{for } 0 \leq t < \zeta^0, \\ \partial, & \text{for } t \geq \zeta^0, \end{cases}$$

and $\mathcal{M}^0$ is the $\sigma$-field generated by $X_t^0$ with a usual augmentation by null sets.

The process $X^0$ is a standard process on $E_0$ and it is convenient to introduce the following related notions. A nearly Borel set $A \subset E_0$ is called $X^0$-invariant if $\mathbf{P}_x(\Omega_A) = 1$ for every $x \in A$, where

$$\Omega_A = \{\omega \in \Omega : X_t^0, X_{t-}^0 \in A \text{ for every } t \in [0, \zeta^0)\}.$$

Then the restriction $X^0|_A$ defined in a natural way is a standard process on $A$. We say a random variable $\gamma$ on $\Omega$ is $X|_A$-measurable if the restriction $\gamma|_{\Omega_A}$ is measurable with respect to the $\sigma$-field $\mathcal{M}^0 \cap A$. The random variable $\gamma$ need not be defined on $\Omega \setminus \Omega_A$ in this case.

A nearly Borel set $N \subset E_0$ is called $X^0$-properly exceptional if $E_0 \setminus N$ is $X^0$-invariant and $m(N) = 0$. Such a set $N$ is then $X^0$-exceptional in the sense of Lemma 2.2(i).

Throughout this section, we shall assume that $X^0$ admits no killings inside $E_0$; that is,

(3.1) $\qquad \mathbf{P}_x(X_{\zeta^0-}^0 \in E_0, \zeta^0 < \infty) = 0 \qquad$ for every $x \in E_0$,

or, equivalently,

(3.2) $\qquad \kappa_0(dx) := \kappa(dx)|_{E_0} + N(x, F)\mu_H(dx)|_{E_0} = 0.$

Here we state a variant of [6], Definition 1.4. We call a random variable $\gamma = \gamma(\omega)$ on $\Omega$ a *terminal random variable* if there exists an $X^0$-properly exceptional set $N \subset E_0$ such that $\gamma$ is $X^0|_{E_0 \setminus N}$-measurable and

$$\mathbf{E}_x[|\gamma|] < \infty \qquad \text{for every } x \in E_0 \setminus N$$

and

$$\gamma(\theta_t(\omega)) = \gamma(\omega) \qquad \text{for every } \omega \in \Omega_{E_0 \setminus N} \text{ and } t < \zeta^0(\omega).$$

We call such a set $N$ an $X^0$-properly exceptional set for the terminal random variable $\gamma$.

For convenience, let us first make an additional assumption that $E^0$ is an open subset of $E$. This additional condition will be removed in Remark 3.2 and Lemma 3.3.

Then the Dirichlet form $(\mathcal{E}^0, \mathcal{F}^0)$ on $L^2(E_0; m)$ is not only transient by Lemma 2.2(iv) but also regular by [20], Theorem 4.4.3. Let $D_k$ be relatively compact open subsets of $E_0 := E \setminus F$ such that $\overline{D}_k \subset D_{k+1}$ and $F_k \uparrow E_0$. Let $L_k$ be the 0-order equilibrium measure of $D_k$ with respect to $(\mathcal{E}^0, \mathcal{F}^0)$, that is,

$$G^0 L_k(x) = e_k(x) := \mathbf{P}_x(\sigma_{D_k} < \infty) \qquad \text{for } x \in E_0.$$

Note that, for any nonnegative measurable function $f$ on $E_0$,

(3.3) $\qquad \displaystyle\int_{E_0} G^0 f(x) L_k(dx) = \int_{E_0} f(x) e_k(x) m(dx).$

This is because for $f \geq 0$ with $\int_{E_0} f(x) G^0 f(x) m(dx) < \infty$, $G^0 f \in \mathcal{F}_e^0$ and

$$\int_{E_0} G^0 f(x) L_k(dx) = \mathcal{E}(e_k, G^0 f) = \int_{E_0} f(x) e_k(x) m(dx).$$

The general case can be proved by approximation in the same way as in the proof for Lemma 2.3.

LEMMA 3.1. *Assume that $E_0$ is an open subset of $E$.*

(i) *For any $X^0$-excessive function $f$ on $E_0$,*

$$\sup_{k \geq 1} \int_{E_0} f(x) L_k(dx) = L(1, f).$$

(ii) *Assume that condition (3.1) holds. For a terminal random variable $\gamma$ with an $X^0$-properly exceptional set $N \subset E_0$, put*

(3.4) $$f(x) = \mathbf{E}_x[\gamma^2] - (\mathbf{E}_x[\gamma])^2, \qquad x \in E_0 \setminus N.$$

*Then $f$ is excessive with respect to $X^0|_{E_0 \setminus N}$.*

PROOF. (i) Choose nonnegative $h_n$ such that $G^0 h_n \uparrow f$. Then, by (3.3),

$$\sup_{k \geq 1} \int_{E_0} f(x) L_k(dx) = \sup_{k \geq 1} \sup_{n \geq 1} \int_{E_0} G^0 h_n(x) L_k(dx)$$

$$= \sup_{n \geq 1} \sup_{k \geq 1} \int_{E_0} h_n(x) e_k(x) m(dx)$$

$$= \sup_{n \geq 1} \int_{E_0} h_n(x) m(dx).$$

On the other hand,

$$L(1, f) = \sup_{\alpha > 0} \sup_{n \geq 1} \alpha (1 - \alpha G_\alpha^0 \mathbb{1}, G^0 h_n)_{L^2(E_0, m)}$$

$$= \sup_{n \geq 1} \sup_{\alpha > 0} (h_n, \alpha G_\alpha^0 \mathbb{1})_{L^2(E_0, m)}$$

$$= \sup_{n \geq 1} \int_{E_0} h_n(x) m(dx).$$

Therefore, $\sup_{k \geq 1} \int_{E_0} f(x) L_k(dx) = L(1, f)$.

(ii) We put $h(x) = \mathbf{E}_x(\gamma), x \in E_0 \setminus N$. For any relatively compact open set $D \subset E_0$, we define

$$\mathbf{H}_D h(x) := \mathbf{E}_x[h(X_{\tau_D}^0)] \qquad \text{for } x \in E_0 \setminus N,$$

where $\tau_D$ denotes the first exit time from $D$ by $X^0$. Since $X^0$ admits no killings inside $E_0$, we have $\tau_D < \tau_0$ $\mathbf{P}_x$-a.s. on $\{\tau_0 < \infty\}$ for every $x \in E_0 \setminus N$ and

$$\mathbf{H}_D h(x) = h(x) \quad \text{and} \quad h(x)^2 \leq \mathbf{H}_D h^2(x), \qquad x \in E_0 \setminus N.$$

Hence, for $x \in E_0 \setminus N$ having $f(x) < \infty$,

$$\mathbf{H}_D f(x) = \mathbf{E}_x[\gamma^2] - \mathbf{H}_D(h^2)(x) \leq \mathbf{E}_x[\gamma^2] - h(x)^2 = f(x).$$

Clearly, $\mathbf{H}_D f(x) \leq f(x)$ holds for those $x \in E_0 \setminus N$ with $f(x) = \infty$. On the other hand, according to [12], Lemma 12.2, both functions $h(x)$ and $\mathbf{E}_x(\gamma^2)$ are excessive and, accordingly, finely continuous with respect to the standard process $X^0|_{E_0 \setminus N}$. Therefore, $f$ is superharmonic in the sense of [12] with respect to this standard process. Using [12], Theorem 12.4, we conclude that $f$ is excessive with respect to $X^0|_{E_0 \setminus N}$, namely,

(3.5) $\quad f(x) \geq 0, \qquad P_t^0 f(x) \uparrow f(x) \qquad \text{as } t \downarrow 0, \text{ for every } x \in E_0 \setminus N.$ $\qquad \square$

REMARK 3.2. We may drop the assumption that $E_0$ is open and can replace the increasing sequence $\{D_k, k \geq 1\}$ of relatively compact open subsets of $E_0$ by an $\mathcal{E}^0$-nest $\{F_k, k \geq 1\}$ consisting of compact subsets of $D_0$ and replace $L_k$ by the 0-order equilibrium measure of $F_k$ in $(\mathcal{E}^0, \mathcal{F}^0)$. The conclusion of Lemma 3.1(i) remains valid.

Under the additional assumption that $E_0$ is open, let $((\mathcal{F}^0)^{\text{ref}}, \mathcal{E}^{\text{ref}})$ be the reflected Dirichlet space of the regular transient Dirichlet space $(\mathcal{F}^0, \mathcal{E}^0)$ as defined in [6], Definitions 1.6 and 3.1. We put

(3.6) $N = \{\gamma : \gamma \text{ is a terminal variable with } L(1, \mathbf{E}_\cdot[\gamma^2] - (\mathbf{E}_\cdot[\gamma])^2) < \infty\}.$

By Lemma 3.1, it holds that

(3.7) $$(\mathcal{F}^0)^{\text{ref}} = \mathcal{F}_e^0 + HN,$$

where

(3.8) $\qquad HN = \{h : h(x) = \mathbf{E}_x[\gamma] \text{ for q.e. } x \in E_0 \text{ with } \gamma \in N\}.$

For $f = f_0 + h \in (\mathcal{F}^0)^{\text{ref}}$, where $f_0 \in \mathcal{F}_e^0$ and $h = \mathbf{E}_\cdot[\gamma]$ with $\gamma \in N$,

(3.9) $\qquad \mathcal{E}^{\text{ref}}(f,f) = \mathcal{E}(f_0, f_0) + \tfrac{1}{2} L(1, \mathbf{E}_\cdot[\gamma^2] - (\mathbf{E}_\cdot[\gamma])^2).$

Due to the next lemma, however, the above definition of the space $((\mathcal{F}^0)^{\text{ref}}, \mathcal{E}^{\text{ref}})$ makes sense without assuming that $E_0$ is open.

LEMMA 3.3. *The second statement* (ii) *of Lemma* 3.1 *remains valid without the additional assumption that $E_0$ is open.*

PROOF. In the preceding proof, we used the fact that $E_0$ is a locally compact separable metric space in order to apply Dynkin's [12], Theorem 12.4 directly. Without assuming that $E_0$ is open, however, $E_0$ is related to such a nice space by a quasi-homeomorphism. Indeed, the Dirichlet form $(\mathcal{E}^0, \mathcal{F}^0)$ on $L^2(E_0; m)$ is quasi-regular by Lemma 2.2(v) and, accordingly, we can apply [8] to find a regular Dirichlet space $(E', m', \mathcal{F}', \mathcal{E}')$ such that $E_0$ and $E'$ are quasi-homeomorphic: there exist an $\mathcal{E}^0$-nest $\{F_n\}$ on $E_0$, an $\mathcal{E}'$-nest $\{F'_n\}$ on $E'$ and a one to one mapping $q$ from $E_{00} = \bigcup_{n=1}^\infty F_n$ to $\bigcup_{n=1}^\infty F'_n$ with the restriction of $q$ on each $F_n$ being homeomorphic to $F'_n$. Further, $m', \mathcal{F}', \mathcal{E}'$ are the image by $q$ of $m, \mathcal{F}^0, \mathcal{E}^0$, respectively.

Take any terminal random variable $\gamma$. In view of Lemma 2.2(iii), we can find an ($X$-)properly exceptional set $N \subset E_0$ including both $E_0 \setminus E_{00}$ and the $X^0$-properly exceptional set for $\gamma$. Then $N$ is an $X^0$-properly exceptional set. Let $E'_1 = q(E_0 \setminus N)$, and

$$X'_t(\omega) = q(X_t(\omega)), \qquad \omega \in \Omega,$$
$$\mathbf{P}'_{x'}(\Lambda) = \mathbf{P}_{q^{-1}x}(\Lambda), \qquad \Lambda \in \mathcal{M}^0, \qquad x' \in E'_1.$$

Then we see that

$$X' = (\Omega, \mathcal{M}^0, X'_t, \zeta^0, \mathbf{P}'_{x'})_{x' \in E'_1}$$

is a standard process on $E'_1$ and $\gamma$ can be regarded as a terminal random variable with respect to $X'$ with the $X'$-properly exceptional set $N' = E' \setminus E'_1$.

Since $E'$ is a locally compact separable metric space, we conclude in exactly the same way as in the proof of Lemma 3.1 that the function

$$g(x') = \mathbf{E}'_{x'}(\gamma^2) - (\mathbf{E}'_{x'}(\gamma))^2, \qquad x' \in E'_1,$$

is $X'$-excessive: for the transition function $\{P'_t\}$ of $X'$,

$$g(x') \geq 0, \qquad P'_t g(x') \uparrow g(x') \qquad \text{as } t \downarrow 0, x' \in E'_1,$$

which implies the desired property (3.5) of the function $f$ because

$$f(x) = g(qx), \qquad P^0_t f(x) = P'_t g(qx), \qquad x \in E_0 \setminus N. \qquad \square$$

We now remove the additional assumption that $E_0$ is an open subset of $E$. Based on the above lemma, we may and do regard (3.6)–(3.9) as definition of the *reflected Dirichlet space* of the Dirichlet space $(\mathcal{F}^0, \mathcal{E}^0)$. See also [25] for a related approach. Let us introduce the function space $\mathcal{H}_F$ by

$$\mathcal{H}_F = \{\mathbf{H}\varphi|_{E_0} : \varphi \text{ is a measurable function on } F \text{ with } \mathbf{H}|\varphi|(x) < \infty$$
$$\text{for q.e. } x \in E_0 \text{ and } L(1, \mathbf{H}(\varphi^2) - (\mathbf{H}\varphi)^2) < \infty\}.$$

THEOREM 3.4. *Assume condition* (3.1) *holds.*

(i) *The following inclusions hold:*

(3.10) $$\{\mathbf{H}u|_{E_0}: u \in \mathcal{F}_e\} \subset \mathcal{H}_F \subset HN \quad \text{and} \quad \mathcal{F}_e|_{E_0} \subset (\mathcal{F}^0)^{\text{ref}}.$$

(ii) *For $u \in \mathcal{F}_e$,*

(3.11) $$\mathcal{E}(u,u) \geq \mathcal{E}^{\text{ref}}(u|_{E_0}, u|_{E_0}).$$

(iii) *Equality holds in* (3.11) *for every $u \in \mathcal{F}_e$ if and only if*

(3.12) $$\mu_{\langle \mathbf{H}u \rangle}(F) = 0 \quad \text{for every } u \in \mathcal{F}_e.$$

PROOF. (i) Since the Dirichlet form $(\mathcal{E}^0, \mathcal{F}^0)$ is quasi-regular on $L^2(E_0; m)$ by Lemma 2.2, there exists an $\mathcal{E}^0$-nest $\{A_k, k \geq 1\}$ consisting of compact subsets of $E_0$. Under the assumption (3.1) that $X^0$ has no killing inside $E_0$, it holds for q.e. $x \in E_0$ that $\mathbf{P}_x$-a.s. on $\{\sigma_{E \setminus A_k} < \infty\}$,

$$\sigma_{E \setminus A_k} < \sigma_F \quad \text{and} \quad X_{\sigma_{E \setminus A_k}} \in E_0 \setminus A_k.$$

On account of (2.6) and the quasi-left continuity of $X$, we have, for q.e. $x \in E_0$,

$$\lim_{k \to \infty} X_{\sigma_{E \setminus A_k}} = X_{\sigma_F}, \quad \mathbf{P}_x\text{-a.s. on } \{\sigma_F < \infty\}$$

and so it is measurable with respect to $\mathcal{M}^0$. Hence, for every $h = \mathbf{H}\varphi|_{E_0} \in \mathcal{H}_F$, $\gamma := \varphi(X_{\tau_0}) = \mathbb{1}_{\{\sigma_F < \infty\}} \varphi(X_{\sigma_F})$ is a terminal variable of $X^0$ with

$$\mathbf{E}_x[\gamma] = \mathbf{H}\varphi(x) \quad \text{and} \quad L(1, \mathbf{E}.[\gamma^2] - (\mathbf{E}.[\gamma])^2) = L(1, \mathbf{H}(\varphi^2) - (\mathbf{H}\varphi)^2).$$

So $h \in HN$, which proves $\mathcal{H}_F \subset HN$.

To prove that $\{\mathbf{H}u|_{E_0} : u \in \mathcal{F}_e\} \subset \mathcal{H}_F$, take any $u \in \mathcal{F}_e$. We see then, by Lemma 4.6.6 of [20], that $\mathbf{H}|u|(x)$ is finite for q.e. $x \in E$. For $n \geq 1$, define $u_n = ((-n) \vee u) \wedge n$, which is in $\mathcal{F}_{e,b}$. Observe that, in view of Theorem 2.11 and its proof, without loss of generality, we may and do assume $m(E_0) < \infty$. Then, by virtue of Lemma 2.3 and the first part of Lemma 2.4,

$$\mu_{\langle \mathbf{H}u_n \rangle}(E_0) = L(1, w_n) \quad \text{where } w_n = \mathbf{H}(u_n^2) - (\mathbf{H}u_n)^2,$$

which particularly implies

$$\tfrac{1}{2}\alpha(1 - \alpha G_\alpha^0 \mathbb{1}, w_n)_{L^2(E_0;m)} \leq \mathcal{E}(\mathbf{H}u_n, \mathbf{H}u_n) \leq \mathcal{E}(u_n, u_n) \leq \mathcal{E}(u,u).$$

We first let $n \to \infty$ using Fatou's lemma and then let $\alpha \to \infty$ to get

$$\tfrac{1}{2} L(1, \mathbf{H}(u^2) - (\mathbf{H}u)^2) \leq \mathcal{E}(u,u) < \infty,$$

proving that $\mathbf{H}u|_{E_0} \in \mathcal{H}_F$.

That $\mathcal{F}_e|_{E_0} \subset (\mathcal{F}^0)^{\text{ref}}$ is now clear from the decomposition of $u \in \mathcal{F}_e$:

(3.13) $$u = u_0 + \mathbf{H}u \quad \text{with } u_0 = u - \mathbf{H}u \in \mathcal{F}_e^0.$$

(ii) We claim

$$\mathcal{E}(u,u) - \mathcal{E}^{\mathrm{ref}}(u|_{E_0}, u|_{E_0}) \tag{3.14}$$
$$= \tfrac{1}{2}\mu_{\langle \mathbf{H}u \rangle}(F) + \tfrac{1}{2}\int_F u(x)^2 \kappa(dx) \qquad \text{for every } u \in \mathcal{F}_e.$$

For any $u \in \mathcal{F}_e$, we decompose it as in (3.13). By Theorem 4.6.5 of [20],

$$\mathcal{E}(u,u) = \mathcal{E}(u_0, u_0) + \mathcal{E}(\mathbf{H}u, \mathbf{H}u).$$

We know from (i) that $u|_{E_0} \in (\mathcal{F}^0)^{\mathrm{ref}}$ and

$$\mathcal{E}^{\mathrm{ref}}(u|_{E_0}, u|_{E_0}) = \mathcal{E}(u_0, u_0) + \tfrac{1}{2}L(1,w).$$

Suppose $u \in \mathcal{F}_{e,b}$. Then by Lemma 2.3 and the first part of Lemma 2.4,

$$\mu_{\langle \mathbf{H}u \rangle}(E_0) = L(1,w) \qquad \text{where } w = \mathbf{H}u^2 - (\mathbf{H}u)^2.$$

Accordingly,

$$\mathcal{E}(\mathbf{H}u, \mathbf{H}u) = \tfrac{1}{2}\mu_{\langle \mathbf{H}u \rangle}(E) + \tfrac{1}{2}\mu^k_{\langle \mathbf{H}u \rangle}(E)$$
$$= \tfrac{1}{2}L(1,w) + \tfrac{1}{2}\mu_{\langle \mathbf{H}u \rangle}(F) + \tfrac{1}{2}\int_E (\mathbf{H}u)^2(x)\kappa(dx),$$

which proves the claim (3.14).

Now for general $u \in \mathcal{F}_e$, since $u_n = ((-n) \vee u) \wedge n$ is a bounded function in $\mathcal{F}_e$ and so

$$(3.15) \quad \mathcal{E}(u_n, u_n) - \mathcal{E}^{\mathrm{ref}}(u_n|_{E_0}, u_n|_{E_0}) = \tfrac{1}{2}\mu_{\langle \mathbf{H}u_n \rangle}(F) + \tfrac{1}{2}\int_F u_n(x)^2 \kappa(dx).$$

But $u_n$ is $\mathcal{E}$-convergent to $u$ and, moreover, $u_n|_{E_0}$ is $\mathcal{E}^{\mathrm{ref}}$-convergent to $u|_{E_0}$ by virtue of [6]. Since

$$|\mu_{\langle \mathbf{H}u_n \rangle}(F)^{1/2} - \mu_{\langle \mathbf{H}u \rangle}(F)^{1/2}| \leq \mu_{\langle \mathbf{H}(u_n - u) \rangle}(F)^{1/2}$$
$$\leq \mu_{\langle \mathbf{H}(u_n - u) \rangle}(E)^{1/2}$$
$$\leq \mathcal{E}(u_n - u, u_n - u)^{1/2},$$

we arrive at (3.14) for $u$ by passing $n \to \infty$ in (3.15). (ii) now follows immediately from (3.15).

(iii) Note that, for $u \in \mathcal{F}_e$,

$$\tfrac{1}{2}\int_F u(x)^2 \kappa(dx) = \tfrac{1}{2}\int_F \mathbf{H}u(x)^2 \kappa(dx) \leq \mu_{\langle \mathbf{H}u \rangle}(F)$$

and so (iii) follows from (3.15). □

REMARK 3.5. (i) Condition (3.1) is needed for
$$\{\mathbf{H}u|_{E_0} : u \in \mathcal{F}_e\} \subset HN \quad \text{and} \quad \mathcal{F}_e|_{E_0} \subset (\mathcal{F}^0)^{\text{ref}}$$
to hold. For example, let $X$ be a spherically symmetric $\alpha$-stable process in $\mathbb{R}^d$ and $D$ is a bounded smooth domain in $\mathbb{R}^d$. Define $\tau_D := \inf\{t > 0 : X_t \notin D\}$. It is well known that
$$\mathbf{P}_x(X_{\tau_D-} \in D \text{ and } \tau_D < \infty) = 1 \quad \text{for every } x \in D,$$
and so condition (3.1) fails with $E = \mathbb{R}^d$ and $E_0 := D$. In this example, it is known that
$$(\mathcal{F}^0)^{\text{ref}} = \mathcal{F}^0 = W_0^{\alpha/2,2}(D) \quad \text{while } \mathcal{F}_e|_{E_0} = W^{\alpha/2,2}(D).$$
Hence, when $\alpha > 1$, $(\mathcal{F}^0)^{\text{ref}}$ is strictly contained in $\mathcal{F}_e|_{E_0}$ and, thus, (3.10) fails.

(ii) Under assumption (3.1), condition (3.12) is equivalent to

(3.16) $$\mu_{\langle u \rangle}(F) = 0 \quad \text{for every } u \in \mathcal{F}_e.$$

Clearly, (3.16) implies (3.12). Assume (3.12) holds. For every $u \in \mathcal{F}_e$, let $u_0 = u - \mathbf{H}u$, which is in $\mathcal{F}_e^0$. It can be shown just as in the proof of Theorem 2.10 that $\mu_{\langle u_0 \rangle}^c(F) = 0$. This together with the assumption (3.1) implies that $\mu_{\langle u_0 \rangle}(F) = 0$. Hence,
$$0 \leq \mu_{\langle u \rangle}(F) = \mu_{\langle u_0 + \mathbf{H}u \rangle}(F) \leq (\sqrt{\mu_{\langle u_0 \rangle}(F)} + \sqrt{\mu_{\langle \mathbf{H}u \rangle}(F)})^2 = 0.$$
This proves that (3.12) implies (3.16) and so these two conditions are equivalent.

In the remainder of this section, we consider the space $\overset{\circ}{S}_F$ defined by (2.30), which is the collection of measures charging no set of zero capacity with quasi-support being equal to $F$ q.e. We fix $\mu \in \overset{\circ}{S}_F$ and let $Y$ be the time changed process of $X$ by the PCAF with Revuz measure $\mu$. The process $Y$ is $\mu$-symmetric and its associated Dirichlet form $(\check{\mathcal{E}}, \check{\mathcal{F}})$ on $L^2(F; \mu)$ and the extended Dirichlet space $\check{\mathcal{F}}_e$ are described by (2.31). In Corollary 2.9 and Theorem 2.11, we have derived an explicit expression of the trace form $\check{\mathcal{E}}$ in terms of the Feller measure $U$ and supplementary Feller measure $V$. We will be concerned with an identification of the trace space $\check{\mathcal{F}}_e$ in terms of $U$ and $V$ under condition (3.1).

To this end, let us introduce the space $\mathcal{G}$ of functions on $F$ with finite Douglas integrals by

(3.17)
$$\mathcal{G} = \Big\{\varphi : \varphi \text{ is measurable with } \mathbf{H}|\varphi|(x) < \infty \text{ for } m\text{-a.e. } x \in E_0 \text{ and}$$
$$\int_{F \times F} (\varphi(\xi) - \varphi(\eta))^2 U(d\xi, d\eta) + \int_F \varphi^2(\xi) V(d\xi) < \infty \Big\}.$$

If we assume both conditions (3.1) and (3.12) hold, then we can well deduce from Corollary 2.9 that

$$\check{\mathcal{F}}_e \subset \mathcal{G} \tag{3.18}$$

and

$$\check{\mathcal{E}}(\varphi,\varphi) = \tfrac{1}{2}\int_{F\times F}(\varphi(\xi)-\varphi(\eta))^2 U(d\xi,d\eta) + \int_F \varphi^2(\xi)V(d\xi) \tag{3.19}$$

for every $\varphi \in \check{\mathcal{F}}_e$.

Two functions $\varphi, \psi$ in $\mathcal{G}$ will be regarded to be equivalent if $\mathbf{H}\varphi = \mathbf{H}\psi$, $m$-a.e. on $E_0$. Owing to the definition of Feller measures $U$ and $V$ given by (2.15) and (2.10), the values of the integrals appearing in (3.17) do not depend on the choice of a representative from each equivalence class.

THEOREM 3.6.  *Assume conditions* (3.1) *and* (3.12) *hold. Assume also that*

$$\{\mathbf{H}u|_{E_0} : u \in \mathcal{F}_e\} = \mathcal{H}_F. \tag{3.20}$$

*Then*

$$\check{\mathcal{F}}_e = \mathcal{G}, \tag{3.21}$$

*and identity* (3.19) *holds for every* $\varphi \in \check{\mathcal{F}}_e$.

PROOF.  It suffices to show that (the equivalence classes of) $\mathcal{G}$ is contained in $\check{\mathcal{F}}_e$. Take any $\varphi \in \mathcal{G}$. Note first that $\mathbf{H}|\varphi|(x)$ is finite q.e. on $E_0$ on account of the $X^0$-excessiveness of $\mathbf{H}|\varphi|$, [20], Lemma 4.1.5 and Lemma 2.2 (iii).

Next put

$$\varphi_n = (-n) \vee (\varphi \wedge n) \quad \text{and} \quad f_n = \mathbf{H}\varphi_n^2 - (\mathbf{H}\varphi_n)^2 \qquad \text{for } n \geq 1.$$

As in the proof of the preceding theorem, we assume, without loss of generality, condition (2.16). We then have, by Lemma 2.3,

$$\alpha(q - \alpha G_\alpha^0 q, f_n)_{L^2(E_0;m)} \leq L(q, f_n) = L(f_n, q) \leq L(\mathbf{H}\varphi_n^2, q) = \int \varphi_n^2 \, dV$$

and, accordingly,

$$\alpha(1 - \alpha G_\alpha^0 \mathbb{1}, f_n)_{L^2(E_0;m)} \leq \alpha(\mathbf{H}^\alpha \mathbb{1}, f_n)_{L^2(E_0;m)} + \int_F \varphi_n^2 \, dV.$$

We can now use Lemma 2.5 and (2.24) to obtain

$$\alpha(1-\alpha G_\alpha^0 \mathbb{1}, f_n)_{L^2(E_0;m)} \leq \int_{\mathbf{F}\times F}(\varphi_n(\xi)-\varphi_n(\eta))^2 U(d\xi,d\eta) + 2\int_F \varphi_n^2(\xi)V(d\xi).$$

We first let $n \to \infty$ using Fatou's lemma and then we let $\alpha \to \infty$ to get

$$L(1, \mathbf{H}(\varphi^2) - (\mathbf{H}\varphi)^2) \leq \int_{F \times F} (\varphi(\xi) - \varphi(\eta))^2 U(d\xi, d\eta) + 2 \int_F \varphi(x)^2 V(d\xi) < \infty,$$

which means that $\mathbf{H}\varphi \in \mathcal{H}_F$. By the assumption (3.20), $\mathbf{H}\varphi = \mathbf{H}u$ $m$-a.e. on $E_0$ for some $u \in \mathcal{F}_e$, namely, $\varphi$ is equivalent to $u|_F \in \check{\mathcal{F}}_e$, as was to be proved. $\square$

COROLLARY 3.7. *Assume conditions* (3.1) *and* (3.12) *hold. Assume also that*

(3.22) $$\mathcal{F}_e|_{E_0} = (\mathcal{F}^0)^{\mathrm{ref}}.$$

*Then the equality* (3.21) *holds together with the identity* (3.19) *holding for every* $\varphi \in \check{\mathcal{F}}_e$.

PROOF. Condition (3.22) is equivalent to $\{\mathbf{H}u|_{E_0} : u \in \mathcal{F}_e\} = HN$, which implies (3.20) because of the inclusion (3.10) and the decomposition (3.13). $\square$

REMARK 3.8. In Sections 2 and 3 we have fixed a general quasi-closed subset $F$ of $E$ satisfying condition (2.1) and considered the space $\overset{\circ}{S}_F$ defined by (2.30), the collection of measures charging no set of zero capacity with quasi-support being equal to $F$ q.e. We have obtained a representation for the trace Dirichlet form $(\check{\mathcal{E}}, \check{\mathcal{F}})$ and for the extended trace Dirichlet space $(\check{\mathcal{E}}, \check{\mathcal{F}}_e)$. In particular, we see that $(\check{\mathcal{E}}, \check{\mathcal{F}}_e)$ is independent of the choice of the measure $\mu \in \overset{\circ}{S}_F$.

Take any $\mu \in \overset{\circ}{S}_F$ and denote by $\widetilde{F}_\mu$ the topological support of $\mu$. Although

$$\mu(\widetilde{F}_\mu \setminus F) = 0, \qquad \mathrm{Cap}(F \setminus \widetilde{F}_\mu) = 0,$$

and the trace Dirichlet form $(\check{\mathcal{E}}, \check{\mathcal{F}})$ can be regarded as a regular Dirichlet form on $L^2(\widetilde{F}_\mu; \mu)$, it may happen that $\mathrm{Cap}(\widetilde{F}_\mu \setminus F) > 0$ (see Example 5.1.2 of [20] for such an example). However, by Theorem 6.2.1(iv) of [20], the set $\widetilde{F}_\mu \setminus F$ always has zero $\check{\mathcal{E}}_1$-capacity with respect to the Dirichlet form $(\check{\mathcal{E}}, \check{\mathcal{F}})$.

Assume now that $F$ is a closed subset of $E$ rather than a general quasi-closed set. Then,

$$\widetilde{F}_\mu \subset F, \qquad \mathrm{Cap}(F \setminus \widetilde{F}_\mu) = 0 \qquad \text{for every } \mu \in \overset{\circ}{S}_F,$$

because the quasi-support of any $\mu \in \overset{\circ}{S}_F$ is q.e. included in $\widetilde{F}_\mu \subset F$.

If for a closed subset $F$ of $E$,

(3.23) \qquad there exists a measure $\nu \in \overset{\circ}{S}_F$ with $\widetilde{F}_\nu = F$,

then the closed set $F$ enjoys a nicer property, namely,

$$\widetilde{F}_\mu = F \qquad \text{for every } \mu \in \overset{\circ}{S}_F,$$

and, consequently, $(\check{\mathcal{E}}, \check{\mathcal{F}})$ is a regular Dirichlet form on $L^2(F;\mu)$ for every $\mu \in \overset{\circ}{S}_F$. Indeed, for any $\mu \in \overset{\circ}{S}_F$, the set $F \setminus \widetilde{F}_\mu$ is relatively open in $F$ and $\nu$ does not charge this set of zero capacity, thus, this set must be empty.

EXAMPLE 3.9. Let $D$ be a bounded domain in $\mathbb{R}^n$. Consider the bilinear form $\mathbf{D}$ defined by the Dirichlet integral

$$\mathbf{D}(u, v) = \int_D \nabla u \cdot \nabla v \, dx.$$

Let $H_0^1(D)$, $\widehat{H}^1(D)$ be the closure of $C_c^1(D)$, $C_c^1(\mathbb{R}^n)|_D$ in the Sobolev space $H^1(D)$, respectively. Here $H^1(D)$ is the space of $L^2$-integrable functions on $D$ with $L^2$-integrable first-order distributional derivatives, equipped with inner product $\mathbf{D}(\cdot,\cdot) + (\cdot,\cdot)_{L^2(D,dx)}$. The bilinear form $(\frac{1}{2}\mathbf{D}, \widehat{H}^1(D))$ is an irreducible strongly local regular Dirichlet form on $L^2(\overline{D}; \mathbb{1}_D(x)\,dx)$ and admits an associated recurrent diffusion process $X$ on $\overline{D}$. With $E = \overline{D}$, $F = \partial D$, $E_0 = D$ and $m(dx) = \mathbb{1}_D(x)\,dx$, this Dirichlet form satisfies not only (3.2) [and, consequently, $X$ satisfies (3.1)], but also (3.12) because

$$\mu_{\langle \mathbf{H}u \rangle} = \mathbb{1}_D(x)|\nabla u(x)|^2 \, dx, \qquad u \in \widehat{H}^1(D)$$

and $m(\partial D) = 0$. The absorbed process $X^0$ obtained from $X$ by killing upon leaving $D$ coincides with the absorbing Brownian motion whose Dirichlet form is $(\frac{1}{2}\mathbf{D}, H_0^1(D))$. Consequently, $\partial D$ satisfies condition (2.1).

Assume that

$$\widehat{H}^1(D) = H^1(D), \tag{3.24}$$

which is satisfied, for instance, when $D$ has continuous boundary (see Theorem 2 on page 14 of [29]). Then condition (3.22) is also satisfied (cf. [6]) and we may well call $X$ *the reflecting Brownian motion* on $\overline{D}$. Corollary 3.7 then applies and we get the description

$$\check{\mathcal{F}}_e = \Big\{ \varphi : \mathbf{H}|\varphi|(x) < \infty \text{ a.e. on } D$$

$$\text{and } \int_{\partial D \times \partial D} (\varphi(\xi) - \varphi(\eta))^2 U(d\xi, d\eta) < \infty \Big\} \tag{3.25}$$

$$\check{\mathcal{E}}(\varphi, \varphi) = \tfrac{1}{2} \mathbf{D}(\mathbf{H}\varphi, \mathbf{H}\varphi)$$

$$= \tfrac{1}{2} \int_{\partial D \times \partial D} (\varphi(\xi) - \varphi(\eta))^2 U(d\xi, d\eta), \qquad \varphi \in \check{\mathcal{F}}_e, \tag{3.26}$$

of the extended trace Dirichlet space $\check{\mathcal{F}}_e$ and the trace Dirichlet form $\check{\mathcal{E}}$ of the Sobolev space $H^1(D)$ on $L^2(\partial D; \mu)$ for any measure $\mu \in \overset{\circ}{S}_{\partial D}$. Since $X$ is conservative and the Lebesgue measure of $D$ is finite, the function $q(x) = \mathbf{P}_x(\sigma_{\partial D} = \infty) = 0$ for every $x \in D$ and, therefore, the supplementary Feller measure $V$ vanishes.

The trace Dirichlet form $(\check{\mathcal{E}}, \check{\mathcal{F}}_e \cap L^2(\partial D; \mu))$ is, in general, not a regular Dirichlet form on $L^2(\partial D; \mu)$. But, if we impose condition (3.23), then it is regular on $\partial D$ for any $\mu \in \overset{\circ}{S}_{\partial D}$. Let $D$ be a Lipschitz domain. Then condition (3.24) is satisfied and, since by [1] the reflecting Brownian motion $X$ has a Hölder continuous transition density function with respect to $\mathbb{1}_D(x)\,dx$, $X$ can be refined to start from every point in $\overline{D}$. Moreover, the boundary $\partial D$ satisfies condition (3.23). This is because, by a result of Dahlberg [9], each harmonic measure in $D$ is mutually absolutely continuous with respect to the surface measure $\sigma$ on $\partial D$. Clearly, $\partial D$ is the topological support of $\sigma$. On the other hand, the surface measure $\sigma$ is mutually absolutely continuous with respect to the measure $\mu$ constructed in the proof of Lemma 2.8 with $g \equiv 1$ and, therefore, $\sigma \in \overset{\circ}{S}_{\partial D}$. In particular, $\partial D$ is a quasi-support of $\sigma$, which proves that condition (3.23) holds. Consequently, the Dirichlet form $(\check{\mathcal{E}}, \check{\mathcal{F}}_e \cap L^2(\partial D; \sigma))$ is regular in $L^2(\partial D; \sigma)$.

Fix some $x_0 \in D$. For $x \in D$, let $K(x, \xi)$ be the density of the harmonic measure $\mathbf{P}_x(X_{\sigma_{\partial D}} \in d\xi)$ with respect to the base harmonic measure $\nu(d\xi) := \mathbf{P}_{x_0}(X_{\sigma_{\partial D}} \in d\xi)$. The function $K(x, \xi)$ is called the classical Poisson kernel for $\frac{1}{2}\Delta$ in $D$, which is continuous on $D \times \partial D$ and is harmonic in $x \in D$ for each fixed $\xi \in \partial D$. So for any bounded function $\varphi$ on $\partial D$,

$$\mathbf{H}\varphi(x) = \mathbf{E}_x[\varphi(X_{\sigma_{\partial D}})] = \int_{\partial D} K(x, \xi)\varphi(\xi)\nu(d\xi).$$

Hence, the Feller kernel is well defined by

$$U(\xi, \eta) = L(K(\cdot, \xi), K(\cdot, \eta)), \qquad \xi, \eta \in \partial D,$$

and $U(\xi, \eta)\nu(d\xi)\nu(d\eta)$ gives the Feller measure $U(d\xi, d\eta)$ in (3.25) and (3.26). Moreover, by the Harnack inequality for positive harmonic functions in $D$, the condition that $\mathbf{H}|\varphi|(x) < \infty$ for some (and, hence, for all) $x \in D$ is equivalent to $\varphi \in L^1(\partial D; \nu)$. So (3.25) can be rewritten as

$$\check{\mathcal{F}}_e = \left\{ \varphi \in L^1(\partial D; \nu) : \int_{\partial D \times \partial D} (\varphi(\xi) - \varphi(\eta))^2 U(d\xi, d\eta) < \infty \right\}.$$

Dahlberg [9] proved that, for a general bounded Lipschitz domain $D$, there is an $\varepsilon > 0$ such that $f(\xi) := \frac{\nu(d\xi)}{\sigma(d\xi)}$ is locally in $L^{2+\varepsilon}(\partial D, \sigma)$ and showed that this result cannot be improved in general. However, when $D$ is a bounded $C^{1,1}$ domain in $\mathbb{R}^n$, by using the two-sided Green function estimates in $D$, it can be shown (cf. the proof of Theorem 3.14 in [22]) that $f$ is bounded between two positive constants. Hence, when $D$ is a bounded $C^{1,1}$ domain in $\mathbb{R}^n$,

$$\check{\mathcal{F}}_e = \left\{ \varphi \in L^1(\partial D; \sigma) : \int_{\partial D \times \partial D} (\varphi(\xi) - \varphi(\eta))^2 U(d\xi, d\eta) < \infty \right\}.$$

EXAMPLE 3.10. Let $D$ be an open $n$-set in $\mathbb{R}^n$ whose boundary $\partial D$ has locally finite $(n-1)$-dimensional Hausdorff measure. Consider for $1 < \alpha < 2$ the space defined by

$$(3.27) \qquad \mathcal{F} = \left\{ u \in L^2(D; dx) : \int_{D \times D} \frac{(u(x) - u(y))^2}{|x - y|^{n+\alpha}} \, dx \, dy < \infty \right\},$$

$$\mathcal{E}(u,v) = \mathcal{A}(n, -\alpha) \int_{D \times D} \frac{(u(x) - u(y))(v(x) - v(y))}{|x - y|^{n+\alpha}} \, dx \, dy,$$
(3.28)
$$u, v \in \mathcal{F},$$

where $\mathcal{A}(n, -\alpha)$ is a positive universal constant given by (1.1) that is relevant to the symmetric $\alpha$-stable process on $\mathbb{R}^n$. We refer the reader to [5] for the following facts. $(\mathcal{E}, \mathcal{F})$ is a regular irreducible Dirichlet form on $L^2(\overline{D}; \mathbb{1}_D(x) \, dx)$ and the associated Hunt process $X$ on $\overline{D}$ may be called a *reflected $\alpha$-stable process*. It is shown in [7] that $X$ has Hölder continuous transition density functions with respect to the Lebesgue measure $dx$ on $\overline{D}$ and, therefore, $X$ can be refined to start from every point in $\overline{D}$. For this example, $E = \overline{D}$, $F = \partial D$, $E_0 = D$ and $m(dx) = \mathbb{1}_D(x) \, dx$. Note that since $D$ is an open $n$-set, $\partial D$ has zero Lebesgue measure. The process $X^0$ obtained from $X$ by killing upon leaving $D$ is the *censored $\alpha$-stable process* with Dirichlet form $(\mathcal{E}, \mathcal{F}^0)$ on $L^2(D)$, where $\mathcal{F}^0$ is the closure of $C_0^1(D)$ in $\mathcal{F}$ with respect to $\mathcal{E}_1$. Let $(\mathcal{E}^{\mathrm{ref}}, (\mathcal{F}^0)^{\mathrm{ref}})$ be the reflected Dirichlet space of $(\mathcal{E}, \mathcal{F}^0)$. Then

$$(3.29) \quad (\mathcal{F}^0)^{\mathrm{ref}} \cap L^2(D; dx) = \mathcal{F} \quad \text{and} \quad \mathcal{E}^{\mathrm{ref}}(u,v) = \mathcal{E}(u,v) \qquad \text{for } u, v \in \mathcal{F}.$$

The process $X^0$ has no killings inside $D$ and so it satisfies condition (3.1). It is easy to see that the energy measure of any $u \in \mathcal{F}$ is absolutely continuous with respect to $m$ and so $(\mathcal{E}, \mathcal{F})$ satisfies condition (3.12). Suppose $D$ is of finite Lebesgue measure and $1 < \alpha < 2$, then (2.1) holds and we know from (3.29) and [6] that condition (3.22) is satisfied. Therefore, Corollary 3.7 applies to the Dirichlet form $(\mathcal{E}, \mathcal{F})$ given by (3.27) and (3.28). So we have the identification (3.21) and (3.19) of its trace Dirichlet space $(\check{\mathcal{F}}_e, \check{\mathcal{E}})$ on $\partial D$ in terms of the Feller measure $U$. Note that $X$ is recurrent, as $D$ has finite Lebesgue measure. Since $\partial D$ has positive capacity, $q(x) = \mathbf{P}_x(\sigma_{\partial D} = \infty) = 0$ and so the supplementary Feller measure $V$ vanishes.

Let us assume additionally that $D$ is a bounded $C^{1,1}$-domain. By [22], Theorem 3.14, the surface measure $\sigma$ on $\partial D$ is mutually absolutely continuous with respect to the $X^0$-harmonic measure $\mathbf{P}_x(X_{\sigma_{\partial D}} \in d\xi)$ for every $x \in D$. Moreover, for each $x \in D$, the Radon–Nikodym derivative of $\mathbf{P}_x(X_{\sigma_{\partial D}} \in d\xi)$ with respect to $\sigma(d\xi)$ is bounded between two positive constants. Consequently, just like in Example 3.9, $\sigma \in \overset{\circ}{S}_{\partial D}$. In particular, $\partial D$ is a quasi-support of $\sigma$. Since clearly $\partial D$ is the topological support of the surface measure $\sigma$, condition (3.23) holds. Therefore, the trace Dirichlet form is regular on $L^2(\partial D; \mu)$ for any choice of $\mu \in \overset{\circ}{S}_{\partial D}$. Furthermore, the condition that $\mathbf{H}|\varphi|(x) < \infty$ for some (and, hence, for all) $x \in D$ is equivalent to $\varphi \in L^1(\partial D; \sigma)$. So in this case, (3.21) can be expressed as

$$\check{\mathcal{F}}_e = \left\{ \varphi \in L^1(\partial D; \sigma) : \int_{\partial D \times \partial D} (\varphi(\xi) - \varphi(\eta))^2 U(d\xi, d\eta) < \infty \right\}$$

and (3.19) becomes

$$\check{\mathcal{E}}(\varphi, \varphi) = \tfrac{1}{2} \int_{\partial D \times \partial D} (\varphi(\xi) - \varphi(\eta))^2 U(d\xi, d\eta) \qquad \text{for } \varphi \in \check{\mathcal{F}}_e.$$

By a similar argument as in Example 3.9, it can be shown further that there is a density function $U(\xi, \eta)$ such that $U(d\xi, d\eta) = U(\xi, \eta)\sigma(d\xi)\sigma(d\eta)$.

**4. Excursions and Feller measures.** In this section we shall give a probabilistic characterization of the Feller measure in terms of the trace of endpoints of excursion of $X$ leaving $F$. We will use the same notation as in the previous sections. As before, let $(\mathcal{E}, \mathcal{F})$ be an irreducible regular Dirichlet form on $L^2(E; m)$ and $X$ be an associated $m$-symmetric Hunt process on $E$. The subset $F$ of $E$ is assumed to be quasi-closed and to satisfy condition (2.1). Without loss of generality, we may further assume that $F$ is nearly Borel and finely closed.

Fix $\mu \in \overset{\circ}{S}_F$ and let $A^\mu$ be the PCAF of $X$ associated with it. As is stated in Section 2, we may take as $F$ the support of $A^\mu$. Let $F^r$ denote the set of all points that are regular for $F$; that is, $F^r := \{x \in E : \mathbf{P}_x(\sigma_F = 0) = 1\}$, where $\sigma_F$ is the hitting time of $F$ by $X$. Since $F^r \subset F$ and $\mathrm{Cap}(F \setminus F^r) = 0$ in view of [20], Theorem 4.3.1, we may choose a properly exceptional set $N$ for the PCAF $A^\mu$ containing $F \setminus F^r$ so that, replacing $F$ by $F \setminus N$, we have (cf. [20], Lemma 5.1.11)

$$F = F^r \quad \text{and} \quad \mathbf{P}_x(\sigma_F = \inf\{t : A^\mu_t > 0\}) = 1 \qquad \text{for every } x \in E \setminus N.$$
(4.1)

Let $Y = (Y_t, \mathbf{P}_x)_{x \in F}$ be the time-changed process of $X|_{E \setminus N}$ by the inverse $\{\tau_t, t \geq 0\}$ of $A^\mu$, that is, $Y_t = X_{\tau_t}$, where $\tau_t := \inf\{s \geq 0 : A^\mu_s > t\}$. Then according to [30], (65.9), $Y$ is a $\mu$-symmetric right process with state space $F$. The associated Dirichlet form $(\check{\mathcal{E}}, \check{\mathcal{F}})$ on $L^2(F; \mu)$ is specified by (2.31) and is quasi-regular. Since $(\check{\mathcal{E}}, \check{\mathcal{F}})$ is quasi-homeomorphic to a regular Dirichlet form, we see by Theorem A.1(ii) that $\mathbf{P}_x$-a.s. $Y_{t-} := \lim_{s \uparrow t} Y_s$ exists and takes value in $F$ for every $t \in (0, \zeta^Y)$ for $\check{\mathcal{E}}$-q.e. $x \in F$, and, accordingly, for $\mathcal{E}$-q.e. $x \in F$, in view of [20], Lemmas 6.2.5 and 6.2.8, where $\zeta^Y := A_\infty^\mu$ is the lifetime of $Y$.

By virtue of Theorem A.1 in the Appendix, a Lévy system of $Y$ exists and gives a probabilistic characterization for the Beurling–Deny decomposition for the Dirichlet form $(\check{\mathcal{E}}, \check{\mathcal{F}})$ of $Y$. On the other hand, it follows from (2.35) and (2.36) that the jumping measure $\check{J}$ of $Y$ is

$$\check{J} = \tfrac{1}{2} U + J|_{F \times F} \tag{4.2}$$

and the killing measure $\check{\kappa}$ of $Y$ is

$$\check{\kappa}(dx) = V(dx) + \kappa(dx)|_F, \tag{4.3}$$

where $U$, $V$, $J$ and $\kappa$ are the Feller measure and supplementary Feller measure of $F$, the jumping measure and the killing measure of $X$, respectively. In particular, they are independent of the choice of measure $\mu \in \overset{\circ}{S}_F$.

For any $\omega \in \Omega$, we define

$$M(\omega) = \overline{\{t \in [0, \infty) : X_t(\omega) \in F\}}.$$

Clearly, the relatively open set $M(\omega)^c$ in $[0, \infty)$ consists of all excursion intervals away from $F$ of the sample path $\omega$. We denote by $I$ the set of left endpoints of excursion intervals in $M^c$. $M$ is homogeneous, that is, $M \circ \theta_s + s = M$ if $M \subset [s, \infty)$. $I$ is also homogeneous.

For $t > \sigma_F$, we define

$$L(t) := \sup[0, t] \cap M$$

and

$$R(t) := \inf(t, \infty) \cap M = \inf\{s > t : X_s \in F\},$$

with the convention that $\inf \varnothing = \infty$. When $t > \sigma_F$, we call $(L(t), R(t))$ the excursion straddling on $t$. Clearly, $t \mapsto R(t)$ is right continuous and increasing and it is easy to verify that $R(t) = \sigma_F \circ \theta_t + t$, and that, for any $s, t \geq 0$, $R(t) \circ \theta_s + s = R(t + s)$. Due to (4.1), $X_{R(t)} \in F$ on $\{R(t) < \infty\}$. We can also see that, for $t > \sigma_F$, $R(t-) < R(t)$ if and only if $t \in I$ and in this case $t = R(t-) = L(t)$. We shall further verify in the proof of Theorem 4.2 [see (4.7) below] that $\mathbf{P}_x$-a.s. $X_{R(t-)-} \in F$ for every $t > \sigma_F$ with $R(t) < \infty$ for $\mathcal{E}$-q.e. $x \in E$.

For any nonnegative measurable function $\Psi$ on $E_\partial \times E_\partial$, consider a random measure $\Pi(\Psi, \cdot)$ defined by

$$\Pi(\Psi, dt) = \sum_{0 < s : R(s-) < R(s)} \Psi(X_{R(s-)-}, X_{R(s)}) \varepsilon_s(dt), \tag{4.4}$$

where $\varepsilon_s$ is the point mass at $s$. The random measure $\Pi$ may also be written as

$$\Pi(\Psi, dt) = \sum_{0<s:s\in I} \Psi(X_{L(s)-}, X_{R(s)}) \varepsilon_s(dt).$$

Recall that $I$ is the (random) set of left endpoints of excursion intervals in $M^c$.

LEMMA 4.1. *The random measure $\Pi(\Psi, \cdot)$ is homogeneous for any $\Psi \in \mathcal{B}(E_\partial \times E_\partial)^+$.*

PROOF. Since $R(s) \circ \theta_u + u = R(u+s)$, we have $X_{R(s)} \circ \theta_u = X_{R(u+s)}$ and

$$\begin{aligned}
\Pi(\Psi, dt) \circ \theta_u &= \sum_{u<s+u:\, R(u+s-)<R(u+s)} \Psi(X_{R(u+s)-}, X_{R(u+s)}) \varepsilon_s(dt) \\
&= \sum_{u<s:\, R(s-)<R(s)} \Psi(X_{R(s)-}, X_{R(s)}) \varepsilon_s(dt+u) \\
&= \Pi(\Psi, dt+u).
\end{aligned}$$

This proves the lemma. □

The following is the main result in this section, which asserts that the Feller measure $U(dx, dy)$ is the Revuz measure of $\Pi$ with respect to the measure $m$ and is characterized by the end-points of excursions of the process $X$ leaving $F$.

THEOREM 4.2. *For any nonnegative measurable function $\Psi$ on $F \times F$ vanishing along the diagonal that is extended to be zero outside $F \times F$,*

$$\begin{aligned}
(4.5) \quad & \int_{F \times F} \Psi(x,y) U(dx, dy) \\
&= \uparrow \lim_{t \downarrow 0} \frac{1}{t} \mathbf{E}_m \left[ \sum_{0<s\leq t: R(s-)<R(s)<\infty} \Psi(X_{R(s)-}, X_{R(s)}) \right].
\end{aligned}$$

PROOF. First note that, due to Lemma 4.1, the limit on the right-hand side is an increasing limit.

It follows from (2.35) and its relation with the Lévy system of $Y$ (see Theorem A.1 in the Appendix) that

$$2 \int_{F \times F} \Psi(\xi, \eta) \check{J}(d\xi, d\eta) = \uparrow \lim_{t \downarrow 0} \frac{1}{t} \mathbf{E}_\mu \left[ \sum_{0<s\leq t} \Psi(Y_{s-}, Y_s) \right]$$

$$= \uparrow \lim_{\alpha \uparrow \infty} \alpha \mathbf{E}_\mu \left[ \sum_{0<t<\infty} e^{-\alpha t} \Psi(Y_{t-}, Y_t) \right]$$

$$= \uparrow \lim_{\alpha \uparrow \infty} \alpha \mathbf{E}_\mu \left[ \sum_{0<t<\infty} e^{-\alpha t} \Psi(X_{\tau_t-}, X_{\tau_t}) \right].$$

Now we make a change of variable, replacing $t$ with $A_t^\mu$. On account of (4.1),

$$\tau_{A_t^\mu} = \inf\{s : A_s^\mu > A_t^\mu\} = \inf\{s > t : A_{s-t}^\mu \circ \theta_t > 0\} = \sigma_F \circ \theta_t + t$$

and, accordingly,

(4.6) $\qquad \tau_{A_t^\mu} = R(t) \qquad$ for every $t > 0, \mathbf{P}_x$-a.s. for $x \in E \setminus N$.

This particularly means that $t \mapsto A_t^\mu$ is constant on each excursion interval of sample paths leaving $F$. Here we can further notice the following. Write $\tau_t^- := \tau_{t-}$. Whenever $t > \sigma_F$ and $R(t) < \infty$, we have from (4.6) $\tau_{A_t^\mu}^- = R(t-)$ and so

(4.7) $\quad X_{R(t-)-} = X_{(\tau_{A_t^\mu}^-)-} = Y_{A_t^\mu-} \in F, \qquad \mathbf{P}_x$-a.s. for $\mathcal{E}$-q.e. $x \in E$,

because $Y_{t-} \in F, 0 < t < \zeta^Y, \mathbf{P}_x$-a.s. for $\mathcal{E}$-q.e. $x \in F$, according to the observation made in the paragraph preceding (4.2).

Hence,

$$2 \int_{F \times F} \Psi(\xi, \eta) \check{J}(d\xi, d\eta) = \uparrow \lim_{\alpha \uparrow \infty} \alpha \mathbf{E}_\mu \left[ \sum_{t \in M, R(t) < \infty} e^{-\alpha A_t^\mu} \Psi(X_{R(t-)-}, X_{R(t)}) \right]$$

$$= \uparrow \lim_{\alpha \uparrow \infty} \int_F \alpha \mathbf{E}^x[\Sigma_\alpha] \mu(dx),$$

where

$$\Sigma_\alpha := \sum_{t \in M, R(t) < \infty} e^{-\alpha A_t^\mu} \Psi(X_{R(t-)-}, X_{R(t)}).$$

Let us first assume that $m(E) < \infty$. Since $\mu$ is the Revuz measure of the PCAF $A^\mu$ of the $m$-symmetric process $X$, we have from the above, [20], Theorem 5.1.3 and [30], (32.6) that

$$2 \int_{F \times F} \Psi(\xi, \eta) \check{J}(d\xi, d\eta) = \uparrow \lim_{\alpha \uparrow \infty} \alpha \int_F \mathbf{E}_x[\Sigma_\alpha] \mu(dx)$$

(4.8)
$$= \uparrow \lim_{\alpha \uparrow \infty} \alpha \left( \uparrow \lim_{s \downarrow 0} \frac{1}{s} \mathbf{E}_m \left[ \int_0^s \mathbf{E}_{X_u}[\Sigma_\alpha] \, dA_u^\mu \right] \right)$$

$$= \uparrow \lim_{s \downarrow 0} \frac{1}{s} \left( \sup_{\alpha > 0} \alpha \mathbf{E}_m \left[ \int_0^s \Sigma_\alpha \circ \theta_u \, dA_u^\mu \right] \right).$$

Now

$$\alpha \mathbf{E}_m \left[ \int_0^s \Sigma_\alpha \circ \theta_u \, dA_u^\mu \right]$$

$$= \alpha \mathbf{E}_m \left[ \int_0^s \sum_{t \in M \circ \theta_u, R(t+u) < \infty} e^{-\alpha(A_{t+u}^\mu - A_u^\mu)} \Psi(X_{R(t+u)-}, X_{R(t+u)}) \, dA_u^\mu \right]$$

$$= \alpha \mathbf{E}_m \left[ \int_0^s e^{\alpha A_u^\mu} \, dA_u^\mu \sum_{t > u, t \in M, R(t) < \infty} e^{-\alpha A_t^\mu} \Psi(X_{R(t)-}, X_{R(t)}) \right]$$

$$= \mathbf{E}_m \left[ \sum_{t \in M, R(t) < \infty} e^{-\alpha A_t^\mu} \Psi(X_{R(t)-}, X_{R(t)}) \int_0^s I_{\{t>u\}} \, de^{\alpha A_u^\mu} \right]$$

$$= E_m \left[ \sum_{t \in M, R(t) < \infty} e^{-\alpha A_t^\mu} \Psi(X_{R(t)-}, X_{R(t)}) \cdot (e^{\alpha A_{s \wedge t}^\mu} - 1) \right]$$

$$= \mathbf{E}_m \left[ \sum_{t \leq s, t \in M, R(t) < \infty} (1 - e^{-\alpha A_t^\mu}) \Psi(X_{R(t)-}, X_{R(t)}) \right]$$

$$+ \mathbf{E}_m \left[ (e^{\alpha A_s^\mu} - 1) \sum_{s < t, t \in M, R(t) < \infty} e^{-\alpha A_t^\mu} \Psi(X_{R(t)-}, X_{R(t)}) \right].$$

Choose $\Psi$ for which the integral $\int_{F \times F \setminus d} \Psi(x,y) \check{J}(dx, dy)$ is finite. Then by (4.8), $\mathbf{E}_\mu(\Sigma_\alpha) < \infty$ for every $\alpha > 0$.

Since $m(E) < \infty$, $m\mathbf{H}$ is in $\overset{\circ}{S}_F$ (see the proof of Lemma 2.8). Since the measure $\check{J}$ is independent of the choice of $\mu \in \overset{\circ}{S}_F$, we may choose $\mu = m\mathbf{H}$. It follows from the fact $A_{\sigma_F}^\mu = 0$ that

$$\Sigma_\alpha = \sum_{t > \sigma_F, t \in M, R(t) < \infty} e^{-\alpha A_t^\mu} \Psi(X_{R(t)-}, X_{R(t)}) = \Sigma_\alpha \circ \theta_{\sigma_F}.$$

Hence, $\mathbf{E}_m(\Sigma_\alpha) = \mathbf{E}_\mu(\Sigma_\alpha) < \infty$. Then we have

(4.9)
$$\alpha \mathbf{E}_m \left[ \int_0^s \Sigma_\alpha \circ \theta_u \, dA_u^\mu \right]$$
$$= \mathbf{E}_m \left[ \sum_{0 < t \leq s, t \in M, R(t) < \infty} \Psi(X_{R(t)-}, X_{R(t)}) \right]$$
$$+ \mathbf{E}_m \left[ \sum_{0 < s < t, t \in M, R(t) < \infty} e^{-\alpha(A_t - A_s)} \Psi(X_{R(t)-}, X_{R(t)}) \right] - \mathbf{E}_m[\Sigma_\alpha]$$

$$= \mathbf{E}_m\left[\sum_{0<t\leq s, t\in M, R(t)<\infty} \Psi(X_{R(t-)-}, X_{R(t)})\right] + \mathbf{E}_m[\Sigma_\alpha \circ \theta_s] - \mathbf{E}_m[\Sigma_\alpha].$$

From the dominated convergence theorem and the fact that we can insert condition $A_t^\mu > 0$ or, equivalently, condition $t > \sigma_F$ into the summand of $\Sigma_\alpha$, it follows that

$$\mathbf{E}_m[\Sigma_\alpha \circ \theta_s] \leq \mathbf{E}_m[\Sigma_\alpha] \quad \text{and} \quad \lim_{\alpha\to\infty} \mathbf{E}_m[\Sigma_\alpha] = 0.$$

Thus, by (4.8) and (4.9),

$$2\int_{F\times F} \Psi(x,y)\check{J}(dx,dy) = \uparrow\lim_{s\downarrow 0} \frac{1}{s}\mathbf{E}_m\left[\sum_{t\leq s, t\in M, R(t)<\infty} \Psi(X_{R(t-)-}, X_{R(t)})\right].$$

The sum on the right-hand side can be divided into two parts: $t \in I$, where $t = R(t-) < R(t)$, and $t \in M \setminus I$, where $t = R(t-) = R(t)$. Note that it follows from (4.7) that $\mathbf{P}_m$-a.s. $X_{R(t-)-} \in F$ for every $t > \sigma_F$ with $R(t) < \infty$. Therefore, we have

$$\lim_{s\to 0} \frac{1}{s}\mathbf{E}_m\left[\sum_{t\leq s, t\in M\setminus I, R(t)<\infty} \Psi(X_{R(t-)-}, X_{R(t)})\right]$$

$$= \lim_{s\to 0} \frac{1}{s}\mathbf{E}_m\left[\sum_{t\leq s, X_{t-}, X_t \in F} \Psi(X_{t-}, X_t)\right] = 2\int_{F\times F} \Psi(\xi,\eta) J(d\xi, d\eta).$$

Finally, using (4.2), we arrive at (4.5).

Now for a general positive Radon measure $m$ on $E$, let $\{D_n, n \geq 1\}$ be a sequence of relatively compact open sets increasing to $E$. For each $n \geq 1$, let $X^{(n)}$ be the subprocess of $X$ killed upon leaving $D_n$. Then by Theorems 4.4.3 and A.2.10 of [20], $X^{(n)}$ is a symmetric Hunt process on $D_n$ with finite symmetrizing measure $\mathbb{1}_{D_n} \cdot m$ and its associated Dirichlet form is regular on $L^2(D; \mathbb{1}_{D_n} \cdot m)$. So formula (4.5) is applicable to each subprocess $X^{(n)}$. Let $U_n$ be the Feller measure for $F \cap D_n$ under the subprocess $X^{(n)}$. By definition (2.15) for the Feller measure $U_n$, we have for

$U_n(f \otimes g)$

$$= \uparrow\lim_{t\downarrow 0} \frac{1}{t}\int_{D_n\setminus F} \mathbf{E}_x[f(X_{\tau_0}^{(n)})] \cdot \mathbf{E}_x[g(X_{\tau_0}^{(n)}); \tau_0 \leq t] m(dx)$$

$$= \uparrow\lim_{t\downarrow 0} \frac{1}{t}\int_{D_n\setminus F} \mathbf{E}_x[f(X_{\tau_0}); \tau_0 < T_n] \cdot \mathbf{E}_x[g(X_{\tau_0}); \tau_0 \leq t, \tau_0 < T_n] m(dx),$$

where $T_n$ is the first exit time of $D_n$. It is obvious that $T_n$ increases a.s. to $\zeta$. Then the above limit is increasing both in $t \downarrow 0$ and in $n \uparrow +\infty$. Taking $n \uparrow \infty$ and switching limits, it follows that $U_n$ increases to $U$. On the other hand, since $m(D_n) < \infty$, (4.5) holds for $X^{(n)}$, that is,

$$\int_{(F \cap D_n) \times (F \cap D_n)} \Psi(x,y) U_n(dx,dy)$$

$$= \uparrow \lim_{t \downarrow 0} \frac{1}{t} \mathbf{E}_{\mathbb{1}_{D_n} \cdot m} \left[ \sum_{0 < s \leq t \colon R(s-) < R(s) < \infty} \Psi(X^{(n)}_{R(s-)-}, X^{(n)}_{R(s)}) \right]$$

$$= \uparrow \lim_{t \downarrow 0} \frac{1}{t} \mathbf{E}_m \left[ \sum_{0 < s \leq t \colon R(s-) < R(s) < T_n} \Psi(X_{R(s-)-}, X_{R(s)}) \right].$$

Note that since $T_n$ increases to $\zeta$, the last expectation above is increasing in $n$. Again, letting $n \uparrow \infty$ and switching two increasing limits on the right-hand side, we have by the monotone convergence theorem identity (4.5). This completes the proof of the theorem. $\square$

We have a similar result for the supplementary Feller measure $V$. For any nonnegative measurable function $f$ on $F$ that is extended to be zero off $F$,

$$\kappa(f \otimes \mathbb{1}_\partial, dt) := \sum_{0 < s < \infty} \mathbb{1}_{\{R(s-) < R(s) = \infty\}} f(X_{R(s-)-}) \varepsilon_s(dt) = f(X_{\gamma-}) \varepsilon_\gamma(dt)$$

is also a homogeneous random measure, where $\gamma$ denotes the last exit time from $F$ by $X$. Using (4.3), it can be shown that

(4.10)
$$\int_F f(\xi) V(d\xi) = \uparrow \lim_{t \downarrow 0} \frac{1}{t} \mathbf{E}_m \left[ \sum_{s \leq t} f(X_{R(s-)-}) \mathbb{1}_{\{R(s-) < R(s) = \infty\}} \right]$$

$$= \uparrow \lim_{t \downarrow 0} \frac{1}{t} \mathbf{E}_m[f(X_{\gamma-}); \gamma \leq t].$$

The proof is similar and we omit the details here.

## APPENDIX: LÉVY SYSTEM AND BEURLING–DENY DECOMPOSITION FOR SYMMETRIC RIGHT PROCESSES

In this Appendix we show that a Lévy system exists for any symmetric right process associated with a quasi-regular Dirichlet form $(\mathcal{E}, \mathcal{F})$ and that the probabilistic characterizations of the Beurling–Deny decomposition (including the jumping measure and killing measure) of $(\mathcal{E}, \mathcal{F})$ remains true for quasi-regular Dirichlet forms and are independent of the choice of a particular process associated with the Dirichlet form. We derive the result by using a quasi-homeomorphism that relates a quasi-regular Dirichlet form to a regular Dirichlet form.

We start with some basic definitions and facts. Suppose that $(\mathcal{E}, \mathcal{F})$ is an $m$-symmetric quasi-regular Dirichlet form on a Hausdorff topological space

$E$, where $m$ is a $\sigma$-finite measure with full support on $E$. For a closed subset $F$ of $E$, let

$$\mathcal{F}_F = \{f \in \mathcal{F} : f = 0 \ m\text{-a.e. on } E \setminus F\}.$$

The following definitions are from [28]:

(1) An increasing sequence of closed sets $\{F_n\}_{n \geq 1}$ of $E$ is an $\mathcal{E}$-*nest* if and only if $\bigcup_{n \geq 1} \mathcal{F}_{F_n}$ is $\mathcal{E}_1$-dense in $\mathcal{F}$, where $\mathcal{E}_1 = \mathcal{E} + (\cdot, \cdot)_{L^2(X,m)}$.
(2) A subset $N \subset E$ is $\mathcal{E}$-*polar* if and only if there is an $\mathcal{E}$-nest $\{F_n\}_{n \geq 1}$ such that $N \subset E \setminus \bigcap_{n \geq 1} F_n$.
(3) A function $f$ on $E$ is said $\mathcal{E}$-*quasi-continuous* if there is an an $\mathcal{E}$-nest $\{F_n\}_{n \geq 1}$ such that $f|_{F_n}$ is continuous on $F_n$ for each $n \geq 1$.
(4) A statement depending on $x \in A$ is said to hold $\mathcal{E}$-quasi-everywhere ($\mathcal{E}$-q.e. in abbreviation) on $A$ if there is an $\mathcal{E}$-polar set $N \subset A$ such that the statement is true for every $x \in A \setminus N$.

A subset $K$ of $E$ is called $\mathcal{E}$-quasi-open (or, $\mathcal{E}$-quasi-closed) if there is an $\mathcal{E}$-nest $\{F_n, n \geq 1\}$ such that $K \cap F_n$ is a relatively open subset (resp. closed subset) of $F_n$ for every $n \geq 1$. This definition is consistent with the one given in Section 2.1 of [20]. For simplicity, we often say quasi-open (or, quasi-closed) instead of $\mathcal{E}$-quasi-open (resp. $\mathcal{E}$-quasi-closed). Clearly, the notions of quasi-open and quasi-closed, in contrast to open and closed, are invariant under quasi-homeomorphisms (cf. [8]).

For a right process $X$ on a state space $E$ with lifetime $\zeta$, let $\partial$ be a cemetery point added to $E$ and define $E_\partial := E \cup \{\partial\}$. Define

$$X_{\zeta-}(\omega) = \begin{cases} \lim_{t \uparrow \zeta(\omega)} X_t(\omega), & \text{if } \zeta(\omega) < \infty \text{ and the limit exists in } E, \\ \partial, & \text{if either } \zeta(\omega) = \infty \text{ or the} \\ & \text{limit does not exists in } E. \end{cases}$$

A Lévy system for $X$ is a pair $(N, H)$, where $N(x, dy)$ is a kernel on $(E_\partial, \mathcal{B}(E_\partial))$ and $H$ is a PCAF of $X$ with bounded 1-potential such that, for any nonnegative Borel function $f$ on $E \times E_\partial$ that vanishes on the diagonal and is extended to be zero elsewhere,

$$\text{(A.1)} \quad \mathbf{E}_x \left( \sum_{s \leq t} f(X_{s-}, X_s) \right) = \mathbf{E}_x \left( \int_0^t \int_{E_\partial} f(X_s, y) N(X_s, dy) \, dH_s \right)$$

for every $x \in E$ and $t \geq 0$.

THEOREM A.1. *Suppose that $(\mathcal{E}, \mathcal{F})$ is an irreducible quasi-regular Dirichlet form on a Hausdorff topological space $E$. Let $X$ be any symmetric Borel right process associated with $(\mathcal{E}, \mathcal{F})$, whose lifetime is denoted by $\zeta$.*

(i) *(Beurling–Deny decomposition).* There is a $\sigma$-finite measure $J$ on $E \times E \setminus d$ and a smooth measure $\kappa$ on $E$ such that, for every $u, v \in \mathcal{F}$,

$$\mathcal{E}(u,v) = \mathcal{E}^c(u,v) + \int_{E \times E \setminus d} (u(x) - u(y))(v(x) - v(y))J(dx,dy)$$
(A.2)
$$+ \int_E u(x)v(x)\kappa(dx),$$

where $\mathcal{E}^c$ is a symmetric form defined on $\mathcal{F}$ satisfying the following strong local property:

$$\mathcal{E}^c(u,v) = 0$$

for any $u, v \in \mathcal{F}$ with $u$ being constant in a quasi-open neighborhood of the quasi-support of $v \cdot m$. Such a decomposition is unique, which is called the Beurling–Deny decomposition of $(\mathcal{E}, \mathcal{F})$. The measures $J$ and $\kappa$ are called the jumping and killing measures, respectively, of $(\mathcal{E}, \mathcal{F})$.

(ii) There is $\Gamma \subset E$ such that $m(\Gamma) = 0$ and $E \setminus \Gamma$ is $X$-invariant in the sense that

$$\mathbf{P}_x(X_{t-} \text{ exists and takes value in } E \setminus \Gamma \text{ for } t \in (0, \zeta) \cup \{\zeta_i\}$$
$$\text{and } X_t \in E_\partial \setminus \Gamma \text{ for every } t \geq 0) = 1$$

for every $x \in E \setminus \Gamma$, and the symmetric right process $X|_{E \setminus \Gamma}$ with state space $E \setminus \Gamma$ is quasi-left continuous on the random time interval $[0, \zeta)$ and has a Lévy system $(N, H)$. Here $\zeta_i$ denotes the finite inaccessible part of $\zeta$. The properly exceptional set is $\mathcal{E}$-polar. So the pair $(N, H)$ will also be called the Lévy system of $X$.

(iii) Let $\mu_H$ denote the Revuz measure of the PCAF $H$ of $X$. Then

(A.3) $\quad J(dx, dy) = \tfrac{1}{2} N(x, dy) \mu_H(dx) \quad and \quad \kappa(dx) = N(x, \partial) \mu_H(dx).$

(iv) For $u \in \mathcal{F}_e$, there is a MAF $M^u$ of $X$ of finite energy and a CAF $N^u$ of $X$ having zero energy such that

$$u(X_t) - u(X_0) = M_t^u + N^u \qquad \text{for } t \geq 0.$$

Let $M^{u,c}$ and $M^{u,d}$ be the continuous and purely discontinuous martingale parts of $M^u$, respectively. Define

$$M_t^{u,k} := (u(X_{\zeta-})\mathbb{1}_{\{\zeta \leq t\}})^p - u(X_{\zeta-})\mathbb{1}_{\{\zeta \leq t\}},$$

where the superscript $p$ stands for the predictable dual projection, and define $M^{u,j} = M^{u,d} - M^{u,k}$. Let $\langle M^{u,c} \rangle$, $\langle M^{u,j} \rangle$ and $\{M^{u,k}\}$ denote the sharp bracket variational processes of the square integrable martingales $M^{u,c}$, $M^{u,j}$ and $M^{u,k}$, respectively. Their corresponding Revuz measures will be denoted as $\mu^c_{\langle u \rangle}$, $\mu^j_{\langle u \rangle}$ and $\mu^k_{\langle u \rangle}$, respectively. Then for $u \in \mathcal{F}_e$,

$$\mu^c_{\langle u \rangle}(E) = 2\mathcal{E}^c(u,u),$$

$$\mu^j_{\langle u \rangle}(E) = 2\int_{E \times E \setminus d} (u(x) - u(y))^2 J(dx, dy),$$

$$\mu^k_{\langle u \rangle}(E) = \int_E u(x)^2 \kappa(dx).$$

PROOF. The theorem is known to be true when $(\mathcal{E}, \mathcal{F})$ is a regular Dirichlet form on a locally compact separable metric space and $X$ is a Hunt process associated with $(\mathcal{E}, \mathcal{F})$. See Theorem 4.5.2, Theorem 4.6.1, Lemmas 5.3.1–5.3.3 and Theorem A.3.21 in [20]. (See also [2].)

We can then use the quasi-homeomorphism technique from [8] to extend these results to quasi-regular Dirichlet form $(\mathcal{E}, \mathcal{F})$ and its associated symmetric process $X$. This technique has now become more or less standard, at least to experts in the field. Nevertheless, for the reader's convenience, we spell out the details below.

For a quasi-regular Dirichlet form $(\mathcal{E}, \mathcal{F})$, it is proved in [8] that $(\mathcal{E}, \mathcal{F})$ is quasi-homeomorphic to a regular Dirichlet form $(\widehat{\mathcal{E}}, \widehat{\mathcal{F}})$ on a locally compact separable metric space $\widehat{E}$ with symmetrizing measure $\widehat{m}$. That is, there are $\mathcal{E}$-nests $\{F_k, k \geq 1\}$ and $\widehat{\mathcal{E}}$-nest $\{\widehat{F}_k, k \geq 1\}$ consisting of compact sets and a one-to-one map $\phi$ from $\bigcup_{k=1}^\infty F_k$ onto $\bigcup_{k=1}^\infty \widehat{F}_k$ such that:

(a) $\phi$ is a topological homeomorphism between $F_k$ and $\widehat{F}_k$ for every $k \geq 1$,

(b) $\widehat{m}$ and $(\widehat{\mathcal{E}}, \widehat{\mathcal{F}})$ are the images of $m$ and $(\mathcal{E}, \mathcal{F})$ under $\phi$, respectively. That is, $\widehat{m} = m \circ \phi^{-1}$ and $\widehat{\mathcal{F}} = \{u \circ \phi^{-1} : u \in \mathcal{F}\}$ and

$$\widehat{\mathcal{E}}(u \circ \phi^{-1}, v \circ \phi^{-1}) = \mathcal{E}(u, v), \qquad u, v \in \mathcal{F}.$$

Let $\widehat{X}$ be a Hunt process associated with the regular Dirichlet form $(\widehat{\mathcal{E}}, \widehat{\mathcal{F}})$ and let $(\widehat{N}, \widehat{H})$ be its Lévy system. Without loss of generality (cf. [20]), we may and do assume that $\bigcup_{k \geq 1} \widehat{F}_k$ is $\widehat{X}$-invariant. Note that by [20] and [28], $\widehat{E} \setminus \bigcup_{k \geq 1} \widehat{F}_k$ and $E \setminus \bigcup_{k \geq 1} F_k$ are $\mathcal{E}$-polar and $\widehat{E}$-polar, respectively. Because the notion of being quasi-open is invariant under quasi-homeomorphism, conclusion (i) holds immediately for $(\mathcal{E}, \mathcal{F})$ by its quasi-homeomorphism to $(\widehat{\mathcal{E}}, \widehat{\mathcal{F}})$. Let $\widetilde{X} = \phi^{-1}(\widehat{X})$ and $\Gamma := E \setminus \bigcup_{k \geq 1} F_k$. Clearly, $E \setminus \Gamma$ is $\widetilde{X}$-invariant and $m(\Gamma) = 0$. By (a)–(b), $\widetilde{X}$ is a symmetric right process associated with Dirichlet form $(\mathcal{E}, \mathcal{F})$ and it has Lévy system $(N, H)$, where

$$N(x, A) := \widehat{N}(\phi(x), \phi(A)) \quad \text{and} \quad H = \widehat{H}.$$

Clearly, $\widetilde{X}|_{E \setminus \Gamma}$ is quasi-left continuous on the random time interval $[0, \zeta)$ and the left-limit of $\widetilde{X}_t$ exists and takes value in $E \setminus \Gamma$. Moreover, $\widetilde{X}_{\zeta_i-}$ exists and takes value in $E \setminus \Gamma$. This establishes (ii) for $\widetilde{X}$. Since (iii) and (iv) hold for the Hunt process $\widehat{X}$ and its regular Dirichlet form $(\widehat{\mathcal{E}}, \widehat{\mathcal{F}})$, it clearly holds for the process $\widetilde{X}$ through the map $\phi$. So we have established (ii)–(iv) for a particular symmetric right process $\widetilde{X}$ associated with the quasi-regular Dirichlet form $(\mathcal{E}, \mathcal{F})$.

Now suppose that $X$ is any other symmetric right process associated with $(\mathcal{E}, \mathcal{F})$. It is known (see [28]) that $X$ and $\widetilde{X}$ are $m$-equivalent; that is, there is a nearly Borel measurable set $S \subset E$ such that:

(c) $m(E \setminus S) = 0$,

(d) $S$ is both $X$-invariant and $\widetilde{X}$-invariant,

(e) When restricted on $S$, the marginal distribution of $X$ is the same as that of $\widetilde{X}$.

From (d)–(e), when restricted to $S$, the process $X$ has the same law as $\widetilde{X}$. So, when restricted to $S$, they have the same Lévy system and, in particular, give the same jumping and killing measures by using the MAF characterization. But from the very definition, $E \setminus S$ does not contribute to the energy measure of $\mu_{\langle u \rangle}$ of $M^u$ (it is defined as a limit under $\mathbf{E}_m$). This shows that (ii)–(iv) hold for the process $X$. This proves the theorem. $\square$

REMARK A.2. Part (i) of Theorem A.1 has been proved in [10]. The main point of our theorem is that a Lévy system exists for any symmetric right process associated with a quasi-regular Dirichlet form and its probabilistic characterization of the Beurling–Deny decomposition (including the jumping measure and killing measure) is independent of the choice of a particular process associated with the Dirichlet form. This is important in Section 4 when we compute excursion laws.

**Acknowledgments.** We thank the referees for their careful reading of this paper.

Z.-Q. CHEN  
DEPARTMENT OF MATHEMATICS  
UNIVERSITY OF WASHINGTON  
SEATTLE, WASHINGTON 98195  
USA  
E-MAIL: zchen@math.washington.edu

M. FUKUSHIMA  
DEPARTMENT OF MATHEMATICS  
KANSAI UNIVERSITY  
OSAKA  
JAPAN  
E-MAIL: fuku@ipcku.kansai-u.ac.jp

J. YING  
DEPARTMENT OF MATHEMATICS  
FUDAN UNIVERSITY  
SHANGHAI  
CHINA  
E-MAIL: jgying@fudan.edu.cn